\documentclass[11pt,reqno]{amsart}
\usepackage[mathscr]{eucal}

\usepackage{xypic}
\usepackage{amsfonts}
\usepackage{amsmath}
\usepackage{amsthm}
\usepackage{amssymb}
\usepackage{latexsym}
\usepackage{eufrak}
\usepackage{fullpage}
\usepackage{palatino}
\usepackage[bookmarks=false]{hyperref}
\usepackage{xr-hyper}
\usepackage{hyperref}

\newtheorem{theorem}[equation]{Theorem}
\newtheorem{proposition}[equation]{Proposition}
\newtheorem{corollary}[equation]{Corollary}
\newtheorem{definition}[equation]{Definition}

\newtheorem{problem}[equation]{Problem}
\newtheorem{lemma}[equation]{Lemma}

\theoremstyle{remark}

\newtheorem*{remark}{Remark}
\makeatletter
\renewcommand{\subsection}{\@startsection{subsection}{2}{0pt}{-3ex
plus -1ex minus -0.2ex}{-2mm plus -0pt minus
-2pt}{\normalfont\bfseries}} \makeatother

\numberwithin{equation}{subsection}

\newcommand{\dis}{{\displaystyle}}

\newcommand{\beq}{\begin{equation}\label}

\newcommand{\eeq}{\end{equation}}

\newcommand{\iso}{{\,\stackrel{_\sim}{\rightarrow}\,}}
\newcommand{\cd}{\!\cdot\!}
\newcommand{\vi}{${\sf {(i)}}\;$}
\newcommand{\vii}{${\sf {(ii)}}\;$}
\newcommand{\viii}{${\sf {(iii)}}\;$}
\newcommand{\iv}{${\sf {(iv)}}\;$}
\newcommand{\sset}{\subset}

\newcommand{\hA}{{\widehat{A}}}
\newcommand{\Amod}{{\mathtt{Mod}}(A:\U_+)}
\newcommand{\AAA}{{\mathtt{Mod}}(A\td A, \fu)}
\newcommand{\AM}{{\mathtt{Mod}}(A,\bar A^+)}
\newcommand{\Amm}{{\mathtt{Mod}}(\sH_c,\bar A^+)}

\newcommand{\BA}{{\mathbb{A}}}
\newcommand{\fu}{{\mathfrak{U}}}
\newcommand{\ssminus}{\smallsetminus}
\newcommand{\id}{{{\mathtt {Id}}}}

\newcommand{\into}{\,\hookrightarrow\,}
\newcommand{\too}{\,\longrightarrow\,}
\newcommand{\onto}{\,\twoheadrightarrow\,}

\DeclareMathOperator{\ad}{{\mathrm{{ad}}}}
\DeclareMathOperator{\ann}{{\mathrm{Ann}}}
\DeclareMathOperator{\gr}{{\mathrm{{gr}}}}
\DeclareMathOperator{\Supp}{{\mathrm{Supp}}}
\DeclareMathOperator{\Hom}{{\mathrm{Hom}}}
\DeclareMathOperator{\gsup}{{\mathcal{S}\mbox{$upp $}}}
\DeclareMathOperator{\rk}{{\mathrm{rk}}}

\DeclareMathOperator{\Spec}{{\mathrm{Spec}}}
\DeclareMathOperator{\sym}{{\mathtt{{Sym}}}}

\newcommand{\irrep}{{\sf{Irrep}}}

\newcommand{\dd}{{\mathscr{D}}}
\newcommand{\La}{\Lambda}

\newcommand{\Ker}{{\mathtt {Ker}}}

\newcommand{\End}{{\mathtt {End}}}
\newcommand{\td}{\mbox{-}}
\newcommand{\opp}{^{\mathtt{op}}}
\newcommand{\Isot}{{\mathtt {Isotropy}}(\G,V)}
\newcommand{\Isotw}{{\mathtt {Isotropy}}(W,\h\oplus\h^*)}
\newcommand{\e}{{\mathsf{e}}}
\newcommand{\la}{\lambda}
\newcommand{\om}{\omega}
\newcommand{\A}{{\mathcal{A}}}

\newcommand{\wh}{\widehat}
\newcommand{\bh}{\delta}
\newcommand{\sing}{{\mathtt{Sing}}}
\newcommand{\lann}{{\mathtt{LAnn}}}

\newcommand{\fm}{{\mathfrak{m}}}
\newcommand{\mm}{{\mathcal{M}}}
\newcommand{\pp}{{\mathscr P}}
\newcommand{\aal}{{A(\leq0)}}
\newcommand{\aall}{{A(<0)}}

\def\C{{\mathbb{C}}}

\def\Z{{\mathbb{Z}}}

\def\oo{{\mathscr O}}
\def\ou{{{}^{^{_\uparrow\!}}\!\mathscr O}}
\def\od{{\mathscr O_{_{\!^{^\downarrow}}}}}
\def\oa{{{}^{^{_\uparrow\!}}\!\mathscr O}_{_{\!^{^\downarrow}}\!}(A\td A)}

\def\U{{\mathcal U}}
\def\UU{{\mathscr{U}}}

\def\BA{{\bar A}}
\def\O{{\mathcal O}(\sH_c)}
\def\HC{{\mathscr{HC}}}

\def\cZ{{\mathcal{Z}}}
\def\sS{{\mathscr{S}}}
\def\var{{\mathscr{V}}}
\def\cI{{\mathcal{I}}}
\def\G{\Gamma}
\def\sH{{\mathsf{H}}}
\def\br{\{-,-\}}
\def\g{{\mathfrak{g}}}
\def\j{{\mathfrak{J}}}
\def\n{{\mathfrak{n}}}
\def\r{{\mathfrak{r}}}

\def\h{{\mathfrak{h}}}
\begin{document}
\title{On  Primitive Ideals}
\author{Victor Ginzburg}
\address{Department of Mathematics, University of Chicago,  Chicago, IL 
60637, USA.}
\email{ginzburg@math.uchicago.edu}
\maketitle
\begin{abstract}
We extend two well-known results
on primitive ideals in enveloping algebras of semisimple Lie algebras,
 the {\it Irreducibility 
theorem} for associated varieties and {\it Duflo theorem}
on primitive ideals,
to  much wider classes of algebras. Our general version of  the Irreducibility
Theorem says that if $A$ is a positively filtered associative
algebra such that $\gr{A}$ is a commutative Poisson algebra
with  {\it finitely many} symplectic leaves, then the associated
variety of any primitive ideal in $A$ is the closure of
a single connected symplectic leaf. Our  general version of  the Duflo theorem
says that if $A$ is an algebra with a "triangular
structure",
  then  any primitive ideal in $A$ is the
annihilator of a simple highest weight module. 
Applications to symplectic reflection algebras
and Cherednik algebras are discussed.
\end{abstract}
{\small
\tableofcontents
}

\section{Introduction}
\subsection{}
Let $A$ be an associative, not necessarily
commutative, algebra.
Recall that a two-sided ideal $I\subsetneq A$ is called {\it primitive}
if it is the annihilator of a simple left $A$-module.
Primitive ideals often play a role similar to
the one played by
maximal ideals in the case of commutative algebras.
The set of all primitive ideals in $A$, equipped with
Jacobson topology, is called the {\it primitive spectrum}
of $A$. It is a
fundamental invariant of $A$ analogous to the affine scheme
$\Spec A$ in the  commutative case.

A complete understanding of the primitive spectrum
is only available in some "very special"
cases. Our motivation in this paper
comes from one such case where $A=\U\g$ is the universal enveloping 
algebra of a complex semisimple Lie algebra $\g$.
Classification of  primitive ideals in $\U\g$
was one of the central themes of representation theory
during the 70--80's, and it is
nowadays  well understood, cf. [J1], [J2] and references therein. Yet, 
details of the classification and the methods involved
still appear to be quite complicated. 

Our goal is to show that,
quite surprisingly, two very important results on
 primitive ideals in $\U\g$ have natural generalizations
to much wider classes of algebras. The first result is the
{\it Duflo Theorem} [Du], stating that any primitive ideal  in $\U\g$
is the
annihilator of a simple {\it highest weight} module.
We generalize this to arbitrary algebras $A$ that have a
triangular structure, i.e., a pair of "opposite"
subalgebras $A^\pm$ such that $A$ is a finite type
$(A^+\td A^-)$-bimodule, satisfying certain additional 
technical conditions, see \S2.
 Our approach is based
on the notion of Jacquet functor borrowed from representation theory,
and on Gabber's {\it Separation Theorem} [Ga].

The second result is a
generalization of  the {\it Irreducibility Theorem}  saying that 
the associated variety of a  primitive ideal in
$\U\g$ is irreducible,
specifically, it is the
closure of a nilpotent conjugacy class in $\g$. The
latter theorem was first partially proved (by a case-by-case
argument) in [BoBr], and in a more conceptual way in [KT]
and [J2] (independently), using
many earlier deep results due to Joseph, Gabber,
Lusztig, Vogan and others. Our 
argument is an adaptation of a much more
direct proof discovered subsequently by Vogan [Vo],
combined with some recent results by Brown--Gordon [BrGo].

\subsection{} We got involved in these matters
because of our interest in 
primitive ideals in symplectic reflection algebras,
a new important class
of associative algebras introduced recently in [EG]
in connection with various   works
 in combinatorics, completely integrable systems, and
generalized McKay correspondence for orbifolds $\C^n/\G$.
In  \S6, we will describe 
all possible subvarieties in  $\C^n/\G$ that may arise as associated
varieties of primitive ideals in a symplectic reflection algebra.

An important special class of   symplectic reflection algebras
comes from Coxeter groups of finite root systems.
The algebras in question  have been 
called {\it rational 
Cherednik algebras},  see [EG], [DO], [BEG],
since they
appeared as a
"rational" degeneration of  {\it double affine 
Hecke algebras} introduced by Cherednik [Ch].
Thus, the latter may (and should) be  thought of as a deformation
of the former. From this point of view,
representation theory of the
rational Cherednik  algebra  is perhaps
 "more fundamental" than (or at least should be studied before)
 that  of the double  affine 
Hecke algebra in the same sense as 
the representation theory of semisimple Lie algebras
is  "more fundamental"  than that  of the corresponding
quantum groups. 

"Experimental evidence" suggests
that although rational Cherednik  algebras are quite far from
enveloping algebras,
 their representation theory 
shares many features of the representation theory
of semisimple Lie algebras, see [BEG].
Thus, classification of
primitive ideals in  Cherednik algebras
seems to be a 
challenging problem.
The results of the present  paper provide a first step in
this direction. 

In more detail,  let $\sH_c$ be the
 Cherednik algebra  associated to the Weyl group $W$ of
a root system in a vector space $\h$.
There is an
analogue
 of the Bernstein--Gelfand--Gelfand
 category of highest weight modules for $\sH_c$, to
be denoted  $\oo_c$, see [BEG].
The isomorphism classes of simple objects in  $\oo_c$ 
are known to be
parametrized by the set $\irrep(W)$ of irreducible representations
of the group $W$. Given $\sigma \in \irrep(W)$,
let $L_\sigma$  denote  the corresponding simple object of  $\oo_c$,
and let $I_\sigma:=\ann(L_\sigma)\sset \sH_c$ denote its annihilator.
The assignment $\Theta: \sigma\mapsto I_\sigma$ gives a map from
$\irrep(W)$ to the set ${\sf{Prim}}(\sH_c)$ of primitive ideals in $\sH_c$,
and our general version of the Duflo theorem implies that
the map $\Theta$ is {\it surjective}, in particular,
 ${\sf{Prim}}(\sH_c)$ is a finite set. Moreover, our irreducibility
theorem implies that $\var(I_\sigma)$, the associated
variety of $I_\sigma,$ is the closure of $W$-saturation
of a certain symplectic leaf in $\h\times\h^*$.
This leads to
two important open problems.
\begin{problem}[Classification of primitive ideals]\label{prob}
\vi Describe the fibers of the map
$\Theta:\irrep(W)\too
{\sf{Prim}}(\sH_c) $
in a more
direct way.

\vii For each  $\sigma \in \irrep(W)$, describe 
the associated variety of $I_\sigma$, in particular,
for which pairs  $\sigma,\tau \in \irrep(W)$ does one have 
$\var(I_\sigma)=\var(I_\tau)$?
\end{problem}
\begin{remark}
The reader should keep in mind that the Cherednik algebra
$\sH_c$ depends on complex parameter(s) "$c$", which is an
analogue of  "central character" in the Lie algebra case.
The solution to Problem \ref{prob} should depend on the value
of  "$c$" in an essential way.
\end{remark}

As a step towards Problem \ref{prob}(ii),
we will compare, in \S6,
the associated variety of
$L_\sigma$ with that of its  annihilator $I_\sigma$.
Specifically, we show (Theorem \ref{levi})
that the associated variety of
any simple highest weight $\sH_c$-module $L$
equals $\var(\ann L)\;\cap\;(\h\times\{0\}),$
the intersection of the associated
variety of the annihilator of $L$ (in $\sH_c$) with
the  Lagrangian subspace $\h\times\{0\}\sset \h\times\h^*$.\medskip

\noindent
{\bf Acknowledgments.} 
          {\footnotesize{     It is a pleasure to acknowledge my
debt to Ofer Gabber for many
 extremely helpful
remarks, and for a very careful reading of the manuscript.}}

\section{Main results}
\subsection{}
Throughout the paper we work with associative algebras over $\C$,
the field of complex numbers.
We write $\otimes=\otimes_\C$.

Let $A$ be a unital $\C$-algebra  equipped
with
a multiplicative increasing filtration:
$0=A_{-1}\subset A_0\subset A_1\subset\ldots,\,$
such that $A_i\cdot A_j\subset A_{i+j}$, and
$\cup_{j\geq 1}\,A_j=A\,,\,1\in A_0$.
Furthermore, we assume that
$\gr A = \bigoplus_j\,A_j/A_{j-1}$,
the associated graded algebra,
is a finitely generated commutative  algebra without zero-divisors,
in other words $\gr A \simeq\C[X]$ is
 the coordinate ring of a (reduced)
irreducible  affine
algebraic variety $X$. For any (say left) ideal $I\sset A$,
the filtration on $A$ induces one on $I$, so that
$\gr I$ becomes an ideal in $\gr A$.
We write $\var(I)\sset X$ for the zero variety
of the ideal $\gr I$, usually referred to as
the associated variety of $I$.

Let $\ell$ be the maximal integer (or $\infty$, if $A$ is commutative) such that,
for all $i,j\geq 0$ and $a_i\in A_i,\,a_j\in A_j,$ one has  
$a_i\cd a_j-a_i\cd a_j \in A_{i+j-\ell}$. (Note that $\ell\geq 1$
since  $\gr A $ is commutative.)
It is well  known that the assignment   $a_i,a_j \mapsto
[a_i,a_j]\,\text{\sl mod}\;A_{i+j-\ell-1}$ descends
to 
a canonical Poisson bracket on $\gr A$ that
makes
 $X$ a 
Poisson algebraic variety. If  $X$ is smooth, then one may view
$X$ as a complex-analytic manifold equipped with a holomorphic
Poisson structure. For each point $x\in X$ one defines
$\sS_x$, the {\it symplectic leaf} through $x$, to be the
set of points that could be reached from $x$  by going along
Hamiltonian flows. 

If  $X$ is not necessarily smooth,
let $\sing(X)$ denote the singular locus of $X$,
and for any $k\geq 1$ define inductively
$\sing^k(X):=\sing\bigl(\sing^{k-1}(X)\bigr)$.
We get a finite partition $X=\bigsqcup_{k=0}^n\;X^k,$
where the strata $X^k:=\sing^{k-1}(X)\ssminus\sing^k(X)$
are smooth analytic varieties (by definition we put
$X^0= X\ssminus \sing(X)$).
It is known, cf. e.g., [BrGo], [Po],
that each $X^k$ inherits a  Poisson structure,
so for any point $x\in X^k$ there is a well  defined
 symplectic leaf $\sS_x\sset X^k$.
This way one defines symplectic leaves on an
arbitrary Poisson algebraic variety. 

In general, each symplectic leaf is a connected
smooth analytic (but not necessarily algebraic)
subset in $X$. However, if the algebraic variety $X$ consists of
finitely many symplectic leaves only, then it was shown
in [BrGo] that each leaf is a smooth irreducible
locally-closed algebraic subvariety in
$X$, and partition into  symplectic leaves
gives an algebraic stratification of $X$.

Our first main result, to be proved in \S3, reads

\begin{theorem}[Irreducibility theorem]\label{irreducibility}
Assume that the Poisson variety $\Spec(\gr{A})$
has only finitely many symplectic leaves. Then,
for any primitive ideal $I\sset A$,
the variety $\var(I)$ is the closure of 
a {\sl single} (connected) symplectic leaf.
\end{theorem}

In the classical case of a semisimple Lie
algebra $\g$, given a primitive ideal
$I\sset\U\g$, let $\cI:=I\cap \cZ$ be the intersection
of $I$ with  $\cZ$, the center of $\U\g$.
A standard argument based on a version of the Schur Lemma, see e.g. [Di],
implies that $\cI$ is a maximal ideal in $\cZ$.
We let $A:= \U\g/\cI\cdot\U\g$ be equipped with the
filtration induced by the standard increasing filtration on $\U\g$.
The algebraic variety
$\Spec(\gr{A})$ is known to coincide with the 
{\it nilpotent variety} in $\g$. The latter variety
is partitioned into finitely many conjugacy classes.
These turn out to be  exactly the symplectic leaves.
Thus, our theorem becomes  the well  known result first
 proved in [BoBr, Theorem 6.5] for classical simple
Lie algebras, and in [J2, Theorem 3.10] and [KT, Proposition 11]
in general
(in the Lie algebra  case our argument reduces to the proof in [Vo],
and does not give anything new).

\subsection{} To formulate our second  result,
let $A=\oplus_{ m\in\Z}\ A( m)$
be a $\Z$-graded unital associative $\C$-algebra.
We fix  a pair  $A^\pm=\oplus_{ m\in\Z}\ A^\pm( m)  \subset A$
of  graded  unital subalgebras.

The above data is assumed to satisfy  the following three conditions:
\smallskip

\noindent
${\sf{(T0)}}\quad$ \parbox[t]{147mm}{$A$ is 
finitely generated as an $(A^+\td A^-)$-bimodule,
i.e., as an $A^+\otimes (A^-)\opp$-module.}

 \noindent
${\sf{(T1)}}\quad$ \parbox[t]{147mm}{The grading on $A$ is
{\em inner}, i.e., there is an element $\bh\in A$ (to be fixed from now on)
such that, for any $ m\in\Z$ and $a\in A( m)$,
one has $[\bh,a]= m\cdot a$.}

\noindent
${\sf{(T2)}}\quad$ \parbox[t]{147mm}{Each of the two  algebras $A^\pm$
  is finitely-generated and  we have $A^+(m)=0$ for all $m<0$, resp.
$A^-(m)=0$ for all $m>0$.
  Furthermore, $A^+(0)=\C=A^-(0)$.}
\smallskip

Let
$\BA^\pm:=\oplus_{m >0}\ A^\pm(\pm m)$ be
the augmentation ideal of the algebra $A^\pm$.

We say that a $\bh$-action  on a vector space $M$
 is locally-finite if any element $m\in M$ is contained
in a $\bh$-stable finite-dimensional vector subspace.
Given a nonunital algebra $\A$ and an $\A$-module $M$,
we say that $\A$-action on
$M$ is locally-nilpotent if, for any $m\in M,$ there exists a large enough
positive integer $n=n(m)$ such that $a_1\cdot \ldots\cdot a_n\cdot m=0\,,\,
\forall a_1,\ldots,a_n\in\A,$
and similarly for right $\A$-modules.

\begin{definition}\label{oo}
Given $(A, A^\pm,\bh)$  as above, let $\ou$
be a full subcategory of the category of 
left $A$-modules whose objects are  finitely-generated
 $A$-modules $M$ such that
the 
$\BA^-$-action  on $M$ is locally-nilpotent.
\end{definition}

In section 4, we will prove the following generalization of
a classical result of J. Bernstein, S. Gelfand, and I. Gelfand \cite{BGG}.

\begin{theorem}\label{criterion}
Let $(A, A^\pm, \bh)$ be a data satisfying conditions
${\sf{(T0)}}$-${\sf{(T2)}}$. Then, one has

\vi The $\bh$-action on any object $M\in\ou$ is
locally-finite and
all (generalized) eigenspaces of this $\bh$-action  are
finite  dimensional.

\vii The category $\ou$ is abelian and any object  of  $\ou$  has finite
length.

\viii The category $\ou$ has only finitely
many (isomorphism classes of) simple objects.

\iv  The category $\ou$
has enough projective and injective objects.
\end{theorem}

\subsection{} 
We call the data $(A, A^\pm,\bh)$ a  {\it commutative triangular structure}
on $A$ if conditions (T0)--(T2) hold and, moreover,
both subalgebras  $A^\pm$ are commutative.
We will also consider noncommutative subalgebras  $A^\pm$.
In that  case, we will assume further that
the algebra $A$ is   equipped
with 
a multiplicative increasing filtration  
$0=A_{-1}\subset A_0\subset A_1\subset\ldots.\,$
We write $\gr A$ for  an associated graded algebra
(which is not assumed to be commutative, in general),
and endow $A^\pm$ with induced filtrations.

We will say
that the data $(A^\pm,\bh)$
 gives a {\it noncommutative triangular structure} on $A$ if,
in addition to (T0)--(T2),
the following holds: 
\smallskip

\noindent
${\sf{(T0')}}\quad$  \parbox[t]{137mm}{$\gr A$ is a
finitely generated $(\gr A^+\td \gr A^-)$-bimodule.}

\noindent
${\sf{(T1')}}\quad$ \parbox[t]{137mm}{We have 
 $\bh\in A_1$, moreover, the image
of $\bh$ in $A_1/A_0$ is a central element in $\gr A$.}

\noindent
${\sf{(T3)}}\quad$ \parbox[t]{137mm}{The algebras
$\gr A^\pm$ are
commutative.}

\noindent
${\sf{(T4)}}\quad$ \parbox[t]{137mm}{The algebras $\gr A^\pm$
are generated by the subspace $\gr_{0}A^\pm\oplus\gr_{1}A^\pm$,
the corresponding degree $\leq 1$ components.
 Moreover, $\dim(\gr_{i}A^\pm) <\infty$, for any $i\geq 0$. }
\smallskip

Clearly, $(T0')$ implies $(T0)$.
Conditions $(T3)$--$(T4)$ imply that the subspaces
$A_1^\pm$ form finite-dimensional Lie subalgebras in $A$
(with respect to the commutator bracket). We denote them
by  $\n^\pm:=A_1^\pm$. It follows, since $\gr A^\pm$ is
generated by the corresponding degree $\leq 1$ component,
that the algebra $A^\pm$ is generated by the subspace
$A_1^\pm$; hence it is the quotient of the augmentation
ideal in  $\U\n^\pm$ by some
other  $\ad\bh$-stable two-sided ideal
in the enveloping algebra $\U\n^\pm$.
Condition $(T1')$ implies that,
for each $i\geq 1,$ the space $A_i$ is 
$\ad\bh$-stable.
Further, the Lie algebras $\n^\pm$ are
necessarily nilpotent, and also
 $1\not\in \bar A^\pm$ (unless $A=0$), because of $(T2)$.

Recall that a two-sided ideal  $I\subsetneq A$ is called prime,
if for any two-sided ideals $J_1,J_2\sset A$ we have 
$\,J_1\cdot J_2\sset I\enspace\Longrightarrow\enspace
J_1\sset I\enspace\text{or}\enspace J_2\sset I.$
Our final important result, proved in \S5,  reads

\begin{theorem}[Generalized Duflo theorem]\label{duflo}
Let $(A^\pm, \bh)$ be  either a commutative or a noncommutative
triangular structure on
$A$. Then, the following conditions
on a two-sided ideal $I\subsetneq A$ are equivalent:
\vskip 2pt
\vi $I$ is a primitive ideal;

\vii $I$ is a prime ideal;

\viii $I$ is the left annihilator
of a simple object of $\ou$.
\vskip 2pt

\noindent
In addition, the algebra $A$ has finite length
as an $A\otimes A\opp$-module.
\end{theorem}

In the classical case, one considers a semisimple Lie
algebra $\g$ with a fixed triangular decomposition
$\g=\n^+\oplus\h\oplus\n^-$. Given a primitive ideal
$I\sset\U\g$, we define $A:= \U\g/\cI\cdot\U\g,$ 
where $\cI:=I\cap \cZ$ as above. 
Let
$A^\pm$ be the image  in $A$ of $\n^\pm\cdot\U\n^\pm$.
Using the Harish-Chandra isomorphism and
the
Poincar\'e--Birkhoff--Witt theorem, it is not difficult to show 
 that $A$ is a finitely generated  $(A^+\td A^-)$-bimodule,
see e.g. \cite{BL} or \cite{Di}. We equip $A$ with an increasing filtration
induced by the standard  increasing PBW-filtration on $\U\g$.
Let $\bh\in\h$ be the half-sum of positive coroots.
Then the adjoint action of $\bh$ on $\U\g$ is known to 
be diagonal with integral eigenvalues.
Furthermore, the data  $(A^\pm,\bh)$
gives a noncommutative triangular structure on $A$.
The categories $\ou$ and $\od$ are both
equivalent, in this
case, to the standard category $\oo$ (with fixed
central character corresponding to $\cI$) as defined
by Bernstein--Gelfand--Gelfand. Thus, our
 Theorem \ref{duflo} applies and it
provides a new proof (cf. also [J2, n. 2.4] for a similar idea)
of the original Duflo theorem~[Du].

\begin{remark} In the special case $A=\U\g$ above,
our proof  of the  Duflo theorem is 
 perhaps not much simpler than the earlier
proofs, see [Du], [BG], [J1],
because of its heavy dependence on the  difficult 
separation theorem due to
Gabber [Ga]. One
advantage of our approach however is that,
unlike all other proofs, it requires very little
information about the representation theory of 
the algebra~$A$.
\end{remark}

\section{Proof of Theorem \ref{irreducibility}}
\subsection{} 
We 
need to recall some results due to K. Brown and I. Gordon [BrGo].

Let $B$ be a finitely generated Poisson algebra,
with  Poisson  bracket $\br$.
An ideal $J\sset B$ (with respect to the commutative product) 
is called a {\it Poisson ideal} if
it is also an ideal with respect to the Lie bracket
$\br$, i.e., if $\{J,B\}\sset J$. 
For the following
standard result of commutative algebra, see e.g. \cite[Theorem 4.5]{CO},
or \cite[\S3.3.2]{Di}.

\begin{lemma}\label{dixmier}
If $J$ is a  Poisson ideal, then so is its radical
and also all the associated  prime ideals.
\hfill\qed
\end{lemma}

Write $\fm_x$ for the maximal  ideal in $B$ corresponding to
a  closed point $x\in \Spec B$, and recall the notation
$\sS_x\sset \Spec B$ for the symplectic leaf through $x$.
Following \cite{BrGo}, given any ideal $J\sset B$ we let $\pp(J)$
denote the maximal Poisson ideal contained in $J$.
If $J$ is prime, then so is  $\pp(J)$, and one has

\begin{proposition}[\cite{BrGo}, Prop.\,3.7]\label{finite_leaves}
Assume the variety $\Spec B$ consists of finitely many
symplectic leaves. Then, for any closed point $x\in \Spec B$,
the leaf $\sS_x$ coincides with the regular locus of the
zero variety of the ideal $\pp(\fm_x)$.
Thus, $\sS_x$ is a smooth connected locally-closed subvariety
in  $\Spec B$.
\hfill\qed
\end{proposition}

Recall  that the Poisson structure on $B$ restricts
to a (nondegenerate) symplectic structure on each symplectic
leaf. In particular, all symplectic leaves
have even (complex) dimension. Hence, from Lemma 
\ref{dixmier} and Proposition \ref{finite_leaves}
we deduce
\begin{corollary}\label{even_dim}
If the variety $\Spec B$ consists of finitely many
symplectic leaves, then every irreducible
component of the zero variety of any
Poisson ideal $J\sset B$ is the closure  of a
symplectic leaf, in particular, it is even dimensional.
\hfill\qed
\end{corollary}

\subsection{} 
Let $A$ be an associative filtered $\C$-algebra
as in Theorem \ref{irreducibility}. 
We recall 

\begin{lemma}[\cite{BoKr}, Korollar 3.6]\label{GKdim}
If $I\sset J \sset A$ are two-sided ideals, $I\neq J$ and moreover $I$
is prime, then   $\dim\var(J)< \dim\var(I)$.
\hfill\qed
\end{lemma}

Given a finitely generated $A$-module $M$, choose a 
{\it good} filtration (cf. e.g. [Be]) on $M$
and write $\Supp(\gr M)$ for the support
of the corresponding associated graded $\gr A$-module,
 a reduced algebraic subvariety in $\Spec(\gr A)$.
It is well  known (due to Bernstein \cite{Be}) that this subvariety
 is independent of the choice of 
a good filtration on $M$, and it will be denoted
$\gsup M$ below. Note that if
$I\sset A$ is a left ideal, then in the special
case $M=A/I$ we have by definition
$\gsup(A/I)=\var(I)$.

Assume from now on that $\gr{A}$ is a Poisson algebra with finitely many
symplectic leaves.
It is straightforward to
verify that, for any two-sided ideal $I\sset A$,
the associated graded  ideal, $\gr{I},$ is a Poisson ideal in $\gr{A}$.
Thus by Corollary \ref{even_dim}, every irreducible
component of  $\var(I)$, the zero variety  of  $\gr{I}$,
is the closure  of a
symplectic leaf.
\smallskip

\noindent
{\bf Proof of Theorem \ref{irreducibility}.\;}
Let $I$ be a primitive ideal. We have just seen that 
proving the theorem  amounts to showing
that  $\var(I)$ is irreducible. To prove this,
we mimic the argument in ([Vo], \S\S3-4). Put 
$M=A/I$, viewed as an $A$-bimodule, and let $X:= \Spec(\gr{A})$.
The filtration on $A$ induces a filtration
on $M$, and we  view $\gr{M}= 
\gr{A}/\gr{I}$ as a finitely generated $\gr{A}$-module.
Let $\mm$ denote the coherent sheaf on  $X$
corresponding to this  $\gr{A}$-module.
By definition we have  $\Supp\mm=\var(I)$,
and we must show that  this variety is irreducible.
Pick an irreducible component of   $\Supp\mm$ of maximal
dimension, $\dim(\Supp\mm)$. 
Corollary \ref{even_dim} says that there exists a 
 symplectic leaf $\sS\sset X$ such that  this irreducible component
is $\overline{\sS}$,
the Zariski closure of $\sS$.
Furthermore, by Lemma \ref{dixmier}, any imbedded component\footnote{
Here  $\var(I)$ is viewed  as a scheme rather than as a reduced variety.}
of $\var(I)$
 (corresponding to an associated prime of $\gr I$)
 that has nonempty intersection with $\sS$ must
contain the whole of $\sS$.

A key point is that since all  symplectic leaves
are even dimensional we have
\begin{equation}\label{codim2}
\dim\bigl(\overline{\sS}\ssminus\sS\bigr)
\leq \dim\sS -2 = \dim(\Supp\mm) -2\,.
\end{equation}

The idea of the proof of Irreducibility Theorem given in ([Vo], \S4), that we will
follow closely, is to introduce a certain new filtration, $F_\bullet,$ on $M=A/I$
defined "locally" in terms of the stratum $\sS$, and such that
$\gr^{F}M$ gives rise to a coherent sheaf on  $\sS$. Then, since the
boundary of  $\sS$ has codimension greater than or equal to
two by \eqref{codim2},
the global sections of $\gr^{F}M$ over $\sS$ necessarily form
a {\it finitely generated} $\gr A$-module which is, moreover,
automatically supported on $\overline{\sS}$. Hence, the filtration
$F$ is good, and one obtains $\var(I)=\gsup(A/I)=\Supp(\gr^{F}A/I)
=\overline{\sS}$.

We proceed to a more detailed exposition.
First, Proposition \ref{finite_leaves} implies that
$\sS$ is Zariski-open in $\overline{\sS}$, hence also in $\Supp\mm$.
Therefore, $Z:=(\Supp\mm)\ssminus \sS$ is a Zariski closed subset in
 $X$. Thus $X\ssminus Z$ is a Zariski open
subset  in  $X$. 
We  may choose a smaller {\it affine} Zariski open
subset $U\subset X\ssminus Z$ such that $(\Supp\mm)\cap U= \sS$.
Write $\jmath: U \into X$
for the corresponding open imbedding.
Set $\mm_{_U}:=\mm|_U$, a coherent sheaf on $U$.
By construction, the leaf $\sS$ is closed in $U,$ and $\Supp\, \mm_{_U}=
\sS$. Hence, $\Supp(\jmath_*\mm_{_U})=
\overline{\sS}.$ Therefore,
 (\ref{codim2}) yields
$$
\dim\bigl(\Supp(\jmath_*\mm_{_U})\;\cap\;
(X\ssminus U)\bigr) \leq\dim(\Supp\,\mm_{_U}) -2\,.
$$
A standard result of algebraic geometry now says that
$\jmath_*\mm_{_U}$ is a coherent sheaf on $X$. Therefore, we have 
\begin{equation}\label{coh}
\G(U\,,\,\mm_{_U})=\G(X\,,\,\jmath_*\mm_{_U})
\quad\text{is a finitely generated $\gr{A}$-module}.
\end{equation}

We are going to apply microlocal techniques --
pioneered by Gabber [Ga1] --
 equivalent to the argument in 
([Vo], \S3), that has been
 also inspired by Gabber. To this end, observe  that the
variety $X$ has a cone-structure, that is,
a contracting $\C^*$-action induced by the 
(nonnegative) grading on $\gr{A}$. Clearly,
$\var(I)$, hence also $\sS$ and $Z$,  are $\C^*$-stable 
subvarieties in~$X$. 

For any nonzero homogeneous element $f\in \gr{A}$,
let $\UU_{\!_{^f}}$
denote   the complement in
$X$ of the divisor $f=0$. Thus, $\UU_{\!_{^f}}$
is an affine Zariski open $\C^*$-stable subset of $X$. 
Let $\wh{A}_{\UU_{\!_{^f}}}$
 be     the formal micro-localization of
$A$ at $\UU_{\!_{^f}}$,   see e.g.  ([Gi], \S1.3).
Thus, $\wh{A}_{\UU_{\!_{^f}}}$    is a complete $\Z$-filtered ring 
such that $\gr(\wh{A}_{\UU_{\!_{^f}}}) =\C[\UU_{\!_{^f}}],$ 
the ring of regular functions
on the affine
open set $\UU_{\!_{^f}}$. The assignment  
$\UU_{\!_{^f}}\longmapsto   \wh{A}_{\UU_{\!_{^f}}}$ extends to
a sheaf $\UU\mapsto \wh{A}_{\UU}$
in
 the topology of Zariski open cone-subsets
in $X\ssminus\{0\}$. For any
Zariski open cone-subset $\UU$, we have an imbedding
$\gr(\wh{A}_{\UU}) \into \C[\UU],$
which becomes an equality if $\UU$ is affine.
There is
a canonical algebra imbedding $A\into \wh{A}_\UU$,
 which is strictly compatible with the filtrations
and such that the induced map
$\gr A \too \gr(\wh{A}_\UU)$ is identified  with
restriction of functions $\gr A=\C[X]\too\C[\UU].$

Further, given a finitely generated $A$-module $L$,
we let $\wh{L}_{\UU_{\!_{^f}}}=\wh{A}_{\UU_{\!_{^f}}}\otimes_A L$ denote its
formal micro-localization, a finitely generated
$\wh{A}_{\UU_{\!_{^f}}}$-module. The assignment   $\UU_{\!_{^f}}\mapsto \wh{L}_{\UU_{\!_{^f}}}$
extends to  a sheaf $\wh{L}$ in
 the topology of Zariski open cone-subsets
in $X\ssminus\{0\}$. One has a  canonical sheaf isomorphism
$\wh{L}\simeq \wh{A}\otimes_A L$; hence, in particular, an isomorphism $\Gamma(\UU,
\wh{L})
\simeq \Gamma(\UU,
\wh{A})\otimes_A L$, for any {\it affine} open cone-subset
$\UU$. 

We now set $\UU:=U$, the affine Zariski open
subset in $X\ssminus
Z$ that has been chosen earlier, cf. \eqref{coh},
and which we now additionally assume to be a cone-subset.
We
apply the formal micro-localization at $U$ to our
$A$-module $M=A/I$ to get the left
$\wh{A}_U$-module $\wh{M}_U=\wh{A}_U\otimes_A (A/I)=
\wh{A}_U/\wh{A}_U\cd I$. We note that since $M$ has
also a right $A$-action (the ideal $I$ is a
two-sided ideal) the left $\wh{A}_U$-module $\wh{M}_U$
acquires a canonical  $\wh{A}_U$-bimodule structure.
Equivalently, one can get the same bimodule
by applying the formal micro-localization procedure
to $M$ viewed as a left $A\otimes A\opp$-module
(from this point of view, $\wh{M}_U$ is "supported" on the 
diagonal  $U \sset U\times U$,
since $M$, viewed as a $A\otimes A\opp$-module,
is clearly  "supported" on the 
diagonal  $X\sset X\times X$).
Further, we have a canonical
$A$-bimodule map
$$
i:\;M\too \wh{A}_U\otimes_A M=\wh{M}_U\,.
$$
Let $K$ denote the kernel of $i$ and $\overline{M}$ denote
the image of $i$. Both  $K$ and $\overline{M}$
are $A$-bimodules again, in particular,
 the preimage of $K$ under the projection
$A\onto A/I=M$ is clearly a two-sided ideal,
say $J\sset A$. Thus, we have the following diagram
\begin{equation}\label{micro_imbed}
{
\diagram
0 \rto & J/I\ddouble\rto & A/I\ddouble\rto&A/J\ddouble\rto&0\\
      & K\rto & M\rto &\overline{M}&\!\!\!\! \!\!\!\!\! \!\!\!
\stackrel{\imath}{\into}\;\wh{M}_U=\wh{A}_U\otimes_A M\,,
\hphantom{xx}
\enddiagram
}
\end{equation}
where we have factored the map $i$ as a composition
$M\onto \overline{M}\stackrel{\imath}{\into}\wh{M}_U$.
The latter factorization shows that the imbedding
$\imath$ induces an isomorphism
$\wh{A}_U\otimes_A\overline{M}\iso\wh{A}_U\otimes_A M$
 of the corresponding micro-localizations.
The  micro-localization functor being exact,
this yields $\wh{K}_U=\wh{A}_U\otimes_A~K=0$.
Hence $(\gsup K)\cap U=\emptyset$.

We claim that $K=0$. If not,
then the two-sided ideal  $J \sset A$, see (\ref{micro_imbed}),
properly contains $I$. Hence, Lemma \ref{GKdim}
yields $\dim(\gsup A/J) <\dim(\gsup A/I)$.
But since $\overline{\sS}$ is a component
of $\gsup A/I$ of maximal dimension,
we get   $\dim\overline{\sS}>\dim(\gsup A/J)$.
Therefore, the equality of sets
$\gsup A/I= (\gsup J/I)\;\cup\;(\gsup A/J),$
which follows from 
the short exact sequence in the top row of (\ref{micro_imbed}),
cannot be achieved unless $\sS\sset \gsup J/I= 
\gsup K$. This last inclusion contradicts
the conclusion of the previous paragraph,
stating that $(\gsup K)\cap U=\emptyset$.
 Thus, $K=0$ and our claim follows.

To proceed further, observe that
the module $\wh{M}_U$ comes equipped with
a natural increasing filtration $F_\bullet\wh{M}_U$
obtained, essentially, by localizing the standard filtration on $M=A/I$.
Moreover, for the associated graded module corresponding
to the  filtration  $F_\bullet\wh{M}_U$ one has
a canonical  isomorphism $\gr^{\!F}(\wh{M}_U)
\iso \G(U\,,\,\mm_{_U})$ since
 $U$ is  affine. We see that 
$\gr^{\!F}(\wh{M}_U)$ 
is a finitely generated $\gr{A}$-module, due to (\ref{coh}).
Further, the increasing filtration on $\wh{M}_U$
induces a filtration $F_\bullet{M}=
M\;\cap\;F_\bullet\wh{M}_U$ via the  imbedding $i: M\into \wh{M}_U$
(we have shown that $i$ is injective).
Then $\gr^{\!F}M$ is a $\gr{A}$-submodule of 
 $\gr^{\!F}(\wh{M}_U)$, hence is also finitely
generated. Therefore the filtration
$F_\bullet{M}$ is a {\it good} filtration on $M$.
This yields $\gsup M =\Supp\bigl(\gr^{\!F}M\bigr)$,
and we deduce
$$
\sS\;\sset\;\var(I)\;=\;\gsup M\;=\;\Supp\bigl(\gr^{\!F}M\bigr)
\;\sset\;\Supp\bigl(\gr^{\!F}\wh{M}_U\bigr)\;\sset\;
\Supp\bigl( \jmath_*\mm_U\bigr)\; \sset\; \overline{\sS}\,.
$$
It follows that $\var(I)=\overline{\sS},$
and  the theorem is proved.\hfill\qed

\section{Category $\ou$}
\subsection{} 
Fix  a data  $(A, A^\pm,\bh)$  satisfying conditions (T0)-(T2).

Given a 
 $\bh$-action on a vector space
$V$ 
 and $\mu\in\C$, let
  $V_\mu:=
\{v\in V\enspace|\enspace (\bh-\mu\cdot\id)^N\cdot
v=0\enspace\text{for}\enspace
 N\gg 0\}$ denote
 the generalized $\mu$-eigenspace of $\bh$. If the $\bh$-action
on $V$ is locally-finite, then there is a (possibly infinite)
direct sum decomposition $V=\oplus_{\mu\in\C}\;V_\mu.$

Let $M$ be an $A$-module.
We consider the following three properties of  $M$:
\begin{equation}\label{item}
{
\parbox[t]{160mm}{
{\sf 1.}\quad
$M$ {\it is finitely generated over $A^+$.} 
\newline
{\sf 2.}\quad{\it The $\bh$-action on $M$ is locally-finite
and each generalized $\bh$-eigenspace $M_\mu$ is finite 
dimensional.}\newline
{\sf 3.}\quad{\it 
The real parts of the eigenvalues of the $\bh$-action in $M$
are all bounded from below.}
}}
\end{equation}

\noindent 
For $m\geq 1$,  put  $(\bar A^-)^m:=
\bar A^-\cdot\ldots\cdot \bar A^-$ ($m$ factors).

\begin{lemma}\label{fin1} 
Let $M$ be an $A$-module generated by a finite  dimensional
subspace $M_0\sset M$ such that, for some $m\geq 1$, we have 
 $(\bar A^-)^m\cdot M_0=0.$ Then properties \eqref{item}(1)--(3) hold.
\end{lemma}

\begin{proof}  We start with statement \eqref{item}(2).
Property (T0) implies that there exists a graded
finite  dimensional
$\C$-vector subspace $E=\oplus_k\ E(k)\subset A$ such that $A=  A^+\cdot E\cdot
 A^-$.
Separating homogeneous components, we get
\begin{equation}\label{ofer}
A(0) \ =\ 
\sum\nolimits_{ k,\ell\in\Z}\
 A^+(- k-\ell)\cdot E(k)\cdot  A^-(\ell).
\end{equation}

Observe  that since the algebra $\bar A^-$ is negatively graded and finitely
generated,
each homogeneous component of $\bar A^-$ is  finite  dimensional,
furthermore, for any integer $m>0$, 
the space $(\bar A^-)^m$
 has finite codimension in $\bar A^-$. Similar properties hold for
the algebra $A^+$.
It follows that, for 
all $\ell\in\C$ except possibly a finite set, one has  
$\bar  A^-(\ell) \sset (\bar A^-)^m$.  Therefore, \eqref{ofer} implies 
that there exist finite  dimensional subspaces $B_m\subset A(0)$
such that
\begin{equation}\label{key}
A(0) = B_m \,+\, A(0)\cap\bigl( A^+\cdot E\cdot
(\bar A^-)^m\bigr)\quad,\quad
m=1,2,\ldots \,.
\end{equation}
Hence,  since $(\bar A^-)^m\cdot M_0=0$, we get
$A(0)\cdot M_0=B_m\cdot  M_0,$  is a finite-dimensional space.
Since $\bh^n\in A(0)$ for any $n\geq 0$, we conclude
that $A(0)\cdot M_0$ is a finite-dimensional 
vector space and that 
$M$ is generated by this vector space as an $A$-module. It follows
that  the  $\bh$-action on $M$ is locally-finite
since the $\ad\bh$-action on $A$ is semisimple. 

To prove statement \eqref{item}(3), 
choose  $B$, a complementary subspace  to $(\bar A^-)^m$
in $ A^-$
so that $ A^-=B\oplus (\bar A^-)^m$. Then we
have  $M= A\cdot M_0=  A^+\cdot E\cdot\bigl(B\oplus (\bar A^-)^m\bigr)
\cdot M_0= A^+\cdot E\cdot B\cdot M_0.$
Further, the space $E\cdot B\cdot M_0$ is
 finite  dimensional, and the action of
$\bar A^+$ may only increase  the real parts of the eigenvalues of
$\bh$-action in $M$. It follows that
 the real parts of the eigenvalues of $\bh$-action in $M$
are all bounded from below, and \eqref{item}(3) is proved.

Finally, since each homogeneous component of $\bar A^+$ is
finite  dimensional,
we deduce from the equality $M= A^+\cdot (E\cdot B\cdot M_0)$ that
 each $\bh$-eigenspace in $M$ is
finite  dimensional. This completes the proof of \eqref{item}(2).
The same  equality clearly implies that $M$  is
finitely generated over $A^+$. 
\end{proof}

\begin{corollary}\label{fin2} Let $M$ be an $A$-module.
Then,  $M\in\ou$  if and only if $M$ is 
 generated by a finite  dimensional
subspace $M_0\sset M$ such that
 $(\bar A^-)^m\cdot M_0=0$, for some $m\geq 1$. In particular,
for any  $m\geq 1$, the $A$-module $A/A\cdot (\bar A^-)^m$
is an object of category $\ou$.
\end{corollary}
\begin{proof} Assume that $M$ is generated by
a subspace $M_0$, as in the statement of the corollary.
Then Lemma \ref{fin1} implies that statement
\eqref{item}(3) holds. Hence,
since the $\bar A^-$-action decreases  the real parts of the eigenvalues of
the $\bh$-action in $M$, we conclude that the  $\bar A^-$-action
on $M$ is locally-nilpotent. This applies in particular
to the module of the form  $A/A\cdot (\bar A^-)^m$,
for any $m\geq 1$.

Conversely, any object $M\in\ou$ is generated by
a finite  dimensional
subspace $M_0\sset M$. The action of  $\bar A^-$ on $M$
is  locally-nilpotent. Hence, there exists $m\gg 0$ such
that $(\bar A^-)^m\cdot M_0=0$.
\end{proof}

\begin{corollary}\label{fin3} The category $\ou$ is abelian, and
properties \eqref{item}(1)--(3) hold for
any object $M\in\ou$.
\end{corollary}
\begin{proof} The last claim is immediate from 
   Lemma \ref{fin1} and
Corollary \ref{fin2}. In particular, we conclude that any object $M\in\ou$
is finitely generated over $A^+$, due to property \eqref{item}(1).
Now, the algebra $\gr A^+$ is commutative and finitely generated,
hence noetherian. It follows that the algebra $A^+$ is noetherian.
Therefore, any $A$-submodule
of a module  $M\in\ou$ is finitely-generated over $A^+$,
hence also over $A$. Thus, the category $\ou$ is abelian.
\end{proof}

\subsection{} 
We put $\aal:=\oplus_{n\leq0}\ A(n)$, resp.
$A(<0):=\oplus_{n<0}\ A(n)$.
Let 
$I:=A(0)\cap A\cdot A(<0)$.
Since, $A^-\sset \aall$, equation \eqref{key} with $m=1$
shows that $I$ is a subspace of finite codimension
in $A(0)$. Furthermore,
it is clear that $I$ is a {\em  two-sided}  ideal 
of the algebra $A(0)$. Thus, 
$H:=A(0)/I$ is a finite-dimensional algebra.

Given an $H$-module $V$, 
we may (and will) view 
$V$ as an $\aal$-module via
an algebra projection
defined as the following composite
$\aal\onto A(0)=\aal/\aall\onto H=A(0)/I$.
We introduce analogues of Verma modules as follows 
\begin{equation}\label{verma}
\Delta(V)\ :=\ A\otimes_{\aal} V.
\end{equation}

The following lemma is an adaptation of
a classical result of J. Bernstein, S. Gelfand, and
I. ~Gelfand \cite{BGG} to our present setting.

\begin{lemma}\label{isoclass} \vi For any simple 
$H$-module $V$, we have $\Delta(V)\in\ou$. Furthermore,
the module $\Delta(V)$ has a unique simple quotient
$L(V)\in\ou$.

\vii Any simple object of  category $\ou$ is isomorphic
to an object of the form $L(V)$ for some simple $H$-module $V$.
\end{lemma}
\begin{proof} Let $V$ be a simple $H$-module.
Then, $V$ has finite dimension since
$H$ is a  finite-dimensional algebra.
Furthermore, the image of $V$ in $\Delta(V)$
is annihilated by the algebra $\aall$, by definition.
Since $\bar A^-\sset \aall$ we deduce from Corollary \ref{fin2} 
that $\Delta(V)\in\ou$.
The proof that the module $\Delta(V)$ has a unique simple quotient
repeats the classical argument word for word.

To prove part (ii), let $L$ be a  (nonzero) simple object of
category $\ou$ and choose $\mu\in\C$  such that
$L_\mu\neq0$ and such that the real part of $\mu$
is minimal possible. Then, one must have
$\aal L_\mu=0$. Furthermore, since $\bh$ is a central
element of the algebra $A(0)$, we deduce that
$L_\mu$ is an $A(0)$-stable subspace.
Moreover,
the space $L_\mu$ is finite dimensional and it is annihilated
by the action of $I$. Thus, one may view
$L_\mu$ as a finite dimensional $H$-module.
Let $V\subset L_\mu$ be a (nonzero) simple
$H$-submodule. The inclusion $V\into L$ extends
to an $A$-module map $\Delta(V)=A\otimes_{\aal} V\
\onto\ L$. This map is surjective since
the space $V$ generates the $A$-module
$L$, the latter being a simple $A$-module.
We conclude that $L$ is a simple quotient of
$\Delta(V)$.\end{proof}

\begin{lemma}\label{fin5} \vi The number of
(isomorphism classes of) simple objects of category $\ou$
is finite.

\vii Any  object of  category $\ou$ has finite length.
\end{lemma}
\begin{proof}
Part (i) follows from  Lemma 
\ref{isoclass} since $H$, the  algebra
involved in the definition of Verma modules \ref{verma},
is 
finite dimensional, hence it has finitely many nonisomorphic simple $H$-modules.

To prove  (ii), 
 fix $M\in \ou$ and a simple object  $L\in \ou$. We claim that
 there exists an integer
$a({M:L})$ such that the following holds:
\vskip -3mm
\begin{align}\label{text}
&\text{\it For any  descending, resp. ascending, chain of
subobjects: $M=M^0\supset M^1\supset \ldots\supset M^N$,}\nonumber\\ 
&\text{\it  let
$S$ be the set formed by the indices "$i$" such that $M^i/M^{i+1}$
contains  $L$ as }\\
&\text{\it a subquotient. Then, one has }
 \;\#S\leq a({M:L}).\nonumber
\end{align}
To prove the claim, choose  $\mu\in\C$
such that the eigenspace $L_\mu$ is nonzero.
Then, we clearly have
$$\dim M_\mu\, =\, \sum\nolimits_{i=1}^N\,\, \dim(M^i/M^{i+1})_\mu
\,\leq\, \sum\nolimits_{i\in S}\,\, \dim(M^i/M^{i+1})_\mu\,\leq\, \#S\cdot \dim
L_\mu\,.
$$
Since $\dim L_\mu\neq 0$, and $\dim M_\mu <\infty$
by part (i) of the theorem, we deduce
$\#S\leq \dim M_\mu/\dim L_\mu$. Thus, we can take $a({M:L})$ to
be the integral part of
$\,(1+\dim M_\mu/\dim L_\mu).$

Now, let   $\{L_1,\ldots,L_d\}$ be  a complete {\it (finite)} collection
of isomorphism classes of simple objects of $\ou$. To prove that
an object $M\in\ou$ has finite length,
one has to verify both ascending and descending chain conditions
for $M$. Consider, for instance, a descending chain of
subobjects
$M=M^0\supsetneq M^1\supsetneq \ldots\supsetneq M^N$.
We must show that
the length of this chain
is bounded from above by a certain number $a=a(M)$ which is independent 
of the chain.
But this is clear from Claim \ref{text}, since for each $i=1,\ldots,N,$
the nonzero object $M^i/M^{i+1}$ must contain
at least one of the simple objects $L_r\,,\,r=1,\ldots,d,$
as a subquotient. Therefore, $N\leq \sum_{r=1}^d\,a({M:L_r})$,
and we may put $a(M):=\sum_{r=1}^d\,a({M:L_r})$.
\end{proof}

\subsection{} For an algebra $A$ satisfying conditions (T0)-(T2),
one can similarly define a category  $\od$ to
be the category of  finitely-generated
{\em right} $A$-modules $M$ such that
the  $\bar A^+$-action  on $M$ is locally-nilpotent.
Below, will also need to  consider a {\em bi}module setting. Specifically,
let $\oa$ be the category of  finitely-generated
left $A\otimes A\opp$-modules $M$ such that the action
on $M$ of the subalgebra
$\bar A^-\otimes A^++
A^- + (\bar A^+)\opp
\sset A\otimes A\opp$ 
is locally-nilpotent.

The above results
concerning category $\ou$ have counterparts for categories
$\od$ and $\oa$.
In particular, we will make use of the following

\begin{lemma}\label{triv} 
\vi
The action of the element
$\bh\otimes 1 + 1\otimes\bh$ on any object $M\in \oa$
is  locally-finite and each 
 generalized eigenspace of the 
$(\bh\otimes 1 + 1\otimes\bh)$-action in  $M\in \oa$
 is finite-dimensional.

\vii Any  simple object
of  category $\oa$ has the form
$L'\boxtimes L'',$
where $L'$ is a simple object in $\ou$
and $L''$  is a simple object in~$\od$.

\viii Category $\oa$ is abelian and
any   object
of that category  has finite length.
\end{lemma}

\begin{proof} Let   $M\in \oa$. Part (i) follows
by an argument  very similar
to the proof of Lemma \ref{fin1}.

We now prove part (ii). Fix a simple object $L\in \oa$.
Choose the nonzero generalized $\bh$-eigenspace of $L$
corresponding to an eigenvalue with minimal possible real part.
Commutation relations show that the action
of any element of  $\bar A^-\otimes A^++
A^- + (\bar A^+)\opp$
decreases the real part of the  eigenvalue. Using this, one
finds a simple $H\otimes H\opp$-module $V$ inside $L$,
cf.   proof of Corollary \ref{fin3}.
The algebra $H$ is finite dimensional,
so the quotient of this algebra by its Jacobson radical
is a direct sum of matrix algebras. It follows
easily that any simple $H\otimes H\opp$-module
has the form $V=V'\boxtimes V''$, where $V'$
is a simple left, resp. $V''$ is a simple right,
$H$-module. Following the proof  of Corollary \ref{fin3},
we deduce that $L$ is a quotient of
an $A\otimes A\opp$-module of the form
$\Delta(V)\boxtimes \Delta\opp(V'')$,
where $\Delta\opp(V''):=A\opp\otimes_{A\opp(\geq0)}V''$, a `right module' counterpart
of the Verma module \eqref{verma}.

We know
  that
$\Delta(V)\in \ou$ and $\Delta\opp(V'')\in \od$ and,
moreover, these objects have finite length, cf. 
 Lemmas \ref{isoclass}-\ref{fin5} and its
$\od$-analogue.
Therefore both $\Delta(V)$
 and $\Delta\opp(V'')$
have finite Jordan--H\"older
series with simple subquotients, say,
$L'_1,\ldots, L'_p$, and $L''_1,\ldots, L''_q,$ respectively.
It is clear from the Jacobson Density Theorem that,
for any $i\in [1,p]\,,\,j\in [1,q],$
the
$A\otimes A\opp$-module $L'_i\boxtimes L''_j$ is simple.
Hence, the object
 $\Delta(V)\boxtimes \Delta\opp(V'')$  also has a finite Jordan--H\"older
series (as an
$A\otimes A\opp$-module) with simple subquotients of the form
$L'_i\boxtimes L''_j.$
It follows that $L$, being a simple quotient of
$\Delta(V)\boxtimes \Delta\opp(V'')$, must be isomorphic to some
simple subquotient from that Jordan-H\"older
series, that is, to some  $L'_i\boxtimes L''_j$, and
part (ii) is proved.

Thanks to part (ii),  category $\oa$ has finitely many nonisomorphic simple
objects. 
Hence, the argument used in the proof of Lemma \ref{fin5}
shows that any  object of category $\oa$ has finite
length. Hence,  any $A\otimes A\opp$-submodule
of an object of category $\oa$ has finite length.
It follows in particular that $\oa$ is an abelian
category.
\end{proof}

Later on, we will need the following 
\begin{lemma}\label{nc} Let $(A^\pm,\bh)$ be a noncommutative
triangular structure on the algebra $A$
and let $M$ be a finitely-generated 
$A\otimes A\opp$-module equipped with a good filtration.
If  $M\in \oa$, then
$\gr M$, an associated graded module,
 is a finitely-generated $\gr(A^+)\otimes \gr(A^-)\opp$-module.
\end{lemma}
\begin{proof} We recall, see e.g. \cite[Corollary 2.3.19]{CG},
 the following well  known result:
if for {\it some} good filtration
on $M$, an associated graded module $\gr M$ is finitely-generated over
$\gr A^+\otimes \gr(A^-)\opp$, then a similar property 
 holds for {\it any} good filtration
on $M$.
Therefore, it suffices to verify the statement of the
lemma   for a particular
  good filtration on $M$ of our choice.
Furthermore, if  $N$ be a subobject of $M$ equipped
with an induced filtration then the validity of
the lemma for the modules $N$ and $M/N$ implies
 the validity of
the lemma for $M$. Since any object of  category
$\oa$ has finite length, we are  reduced to
proving the lemma for simple objects.
The proof of Lemma \ref{triv}
shows that any such module is a quotient of
an  $A\otimes A\opp$-module of the form
$\Delta(V)\boxtimes \Delta\opp(V'')$. The latter
 is itself a quotient of 
$
K:=(A/A\cdot \bar A^-)\,\boxtimes\, (A\opp/A\opp\cdot (\bar A^+)\opp)
$.

Thus,
 it suffices to prove the lemma for
the  module $K$. Furthermore,
we may (and will)  equip  this module  with
a quotient filtration induced by the natural
tensor product
filtration on the algebra $A\otimes A\opp$.
With this choice of filtration, 
using that $\gr(A/A\cdot \bar A^\pm)$
is a quotient of $\gr A\big/\gr A\cdot\gr\bar  A^\pm$,
we deduce that $\gr K$ is a quotient 
 of 
\[{\mathcal K}:=(\gr A\big/\gr A\cdot\gr\bar  A^-)\boxtimes (\gr
A\opp/\gr A\opp\cdot \gr\bar  A^+).
\]

We claim that ${\mathcal K}$ is
a finitely generated $(\gr A^+)\td(\gr
A^-)$-bimodule. To see this, we use property $(T0')$,
which
says that
 $\gr A$ is a finitely generated $(\gr A^+)\td(\gr
A^-)$-bimodule. 
The property implies that $\gr A\big/\gr A\cdot\gr\bar  A^-$
is  a finitely generated $\gr A^+$-module
and $\gr A\opp/\gr A\opp\cdot \gr\bar  A^+$
is  a finitely generated $\gr A^-$-module.
Our claim clearly follows from this.
\end{proof}

\subsection{} 
Given a (possibly infinite-dimensional) vector space $V$
with  locally
finite $\bh$-action such that all generalized $\bh$-eigenspaces $V_\la$
are finite  dimensional,
let $V^\star:= \bigoplus_\la\;(V_\la)^*\sset \Hom_{_\C}(V,\C)$ denote the direct sum of 
the spaces dual to $V_\la$'s. If $V$ is a left
$A$-module, then  $V^\star$ is an $A$-stable subspace
in the right $A$-module $\Hom_\C(V,\C)$
(= linear dual of ~$V$). The right  $A$-module $V^\star$ is
called the  {\it restricted dual} of $V$.
There is
 a canonical  $A$-module isomorphism $V\iso (V^\star)^\star$.

\begin{lemma}\label{dual}
Restricted duality functor
$M\mapsto M^\star$
induces anti-equivalences: $\ou \,
\stackrel{_\sim}{\longleftrightarrow}\, \od.$
\end{lemma}
\begin{proof} 
It is clear that
the real parts of the eigenvalues of the $\bh$-action in $M$
are bounded from below if and only if 
the real parts of the eigenvalues of the $\bh$-action in $M^\star$
are bounded from above. Thus, it follows
by Corollary \ref{fin2} and Lemma \ref{fin1} that   $M\in \ou$
 if and only if $M^\star\in\od$.
\end{proof}

\noindent
{\bf Proof of Theorem \ref{criterion}.\;} 
Part (i) of the
theorem 
is immediate from
Corollary \ref{fin2} and Corollary \ref{fin3}.
Parts (ii) and (iii) follow from Corollary  \ref{fin5}.

To prove part (iv), let $L_1,\ldots, L_d$
be a complete collection of representatives of simple objects
of category $\ou$.  For each
$r=1,\ldots,d,$
there is a real number $\lambda_r\in\mathbb{R}$ such that the real
parts of $\bh$-eigenvalues in $L_r$ are all $\geq \lambda_r$,
by Lemma \ref{fin1}.
It follows  that, for   {\it any} object $M\in\ou$, the real
parts of $\bh$-eigenvalues in $M$
 are all $\geq {\mathtt{min}}(\lambda_1,\ldots,\lambda_d).$
Thus, since the $\bar A^-$-action decreases  the real
parts of $\bh$-eigenvalues by at least $1$, we observe
that, for each $\mu\in \C$ there exists an integer $n(\mu)>0$ such
that we have 
\begin{equation}\label{nmu}
(\bar A^-)^{n(\mu)}\cdot M_\mu=0,\quad\forall M\in \ou\,.
\end{equation}

Given  $\mu\in \C$, for each $k=1,2,\ldots,$ we put
\[
J(\mu,k):= A\cdot (\bh-\mu)^k + A\cdot(\bar A^-)^{n(\mu)} \,.
\]
Clearly, $J(\mu,1)\supseteq J(\mu,2)\supseteq \ldots, $ is a
 descending chain of left ideals in $A$.
We claim that this chain stabilizes, i.e., there exists $\ell\gg 0$
such that we have  $J(\mu,\ell)=J(\mu,\ell+1)=\ldots.$ 
To prove the claim, observe first
that the left $A$-module
$A/A\cdot(\bar A^-)^{n(\mu)}$ is an object of $\ou$,
by Corollary \ref{fin2}. Hence, this object has finite length,
by part (ii) of the theorem. 
Now, for
 $k=1,2,\ldots,$ let $N^k\subset A/A\cdot(\bar A^-)^{n(\mu)}$ be 
an $A$-submodule generated by the coset $(\bh-\mu)^k\;{\sf
{mod}}\,\bigl(A\cdot(\bar A^-)^{n(\mu)}\bigr).$
Clearly $N^1\supseteq N^2\supseteq \ldots.$ We conclude that this
chain of submodules must stabilize; hence there exists $\ell\gg 0$
such that the canonical  projections $
\bigl(A/A\cdot(\bar A^-)^{n(\mu)}\bigr)\big/N^\ell\twoheadleftarrow
\bigl(A/A\cdot(\bar A^-)^{n(\mu)}\bigr)\big/N^{\ell+1}\twoheadleftarrow\ldots$
are all bijections. But these projections are nothing but the
 canonical  projections $A/J(\mu,\ell)\twoheadleftarrow
A/J(\mu,\ell+1)\twoheadleftarrow\ldots,$
and our stabilization claim follows.

We put $P(\mu):= A/J(\mu,\ell)$. The $A$-module $P(\mu)$ is
a quotient of $A/A\cdot(\bar A^-)^{n(\mu)}$, hence,  an object of $\ou$.
We claim\footnote{
This observation has been used earlier, see [Gu].} that 
$P(\mu)$ is a projective in $\ou$. To see this,
let $M\in \ou$ and fix $v\in M_\mu$. Then $(\bar A^-)^{n(\mu)}\cdot v=0$,
due to \eqref{nmu}.
Furthermore, there exists $k$ large enough such that
$(\bh-\mu)^k\cdot v=0$. Therefore, $J(\mu,k)\cdot v=0$,
hence $J(\mu,\ell)\cdot v=0$, by the  stabilization claim proved above. Therefore,
the assignment $a\mapsto a(v)$
descends to a well  defined  $A$-module map $\varphi_v: P(\mu)=A/J(\mu,\ell)
\to M$.
Thus, we get a canonical isomorphism
\begin{equation}\label{functor}
 M_\mu \iso \Hom_A\bigl(P(\mu), M),\quad v\longmapsto
\varphi_v\,.
\end{equation}
We conclude that the functor: $\ou\too \text{\it Vector
spaces},\;
M\mapsto M_\mu$, which is clearly exact,
 is represented by the object $P(\mu)$. Hence,
$P(\mu)$ is projective.

To complete the proof, let $L\in\ou$ be a simple object.
Find $\mu\in \C$ such that $L_\mu\neq 0$. Then 
\eqref{functor} implies that there is a nonzero $A$-module
map $P(\mu)\to L$. This map is surjective since $L$ is simple.
Thus, we have shown that any 
 simple object in $\ou$ is isomorphic to a quotient of 
$P(\mu)$, for an appropriate $\mu\in\C$.
 It follows that the category $\ou$ has enough projectives.
Similarly, one proves that  the category $\od$ has enough projectives
(using Lemma \ref{dual}). But then the equivalence of
Lemma \ref{dual} implies that  category $\ou$ has enough injectives
as well. The theorem is proved.
\qed

\section{Jacquet functor and Proof of Theorem \ref{duflo}}
\subsection{} 
We review the construction of the Jacquet functor following 
[Ca] and [CO].

Let $A$ be an associative algebra and 
${\mathfrak  r}\sset A$ a 
finite-dimensional Lie subalgebra (with respect to the
commutator bracket). Assume further that ${\mathfrak
r}$ has the form of a semidirect product
  ${\n} \rtimes \C \cd \delta$
such that\vskip 3pt

\noindent
$\vi\quad$ The adjoint  action of ${\n}$ on $A$
 is locally-nilpotent and

\noindent
$\vii\quad$
The $\ad\delta$-action on $A$ is locally-finite, and
all eigenvalues of $\ad\delta$-action on $\n$ are in $\Z_+$.
\vskip 3pt

Note that the second condition implies
that $\n$ is a nilpotent Lie algebra.
Let $\U_+$ denote the associative subalgebra  (without unit) in $A$
generated by $\n$. Thus, $\U_+$ is a quotient
of  the augmentation ideal in the enveloping
algebra $\U\n$.

For each $k\geq 1$, let $\U_+^k$ be the $k$-th power of $\U_+$.
These powers form a descending
chain of subalgebras $\U_+\supset \U_+^2 \supset\ldots,$ such that
$\cap_{k\geq 1}\,\U_+^k=0$. Note that condition (i) above implies
that, for any $a\in A$, there exists $n=n(a)\in\Z$ such that
\begin{equation}\label{continuous}
a\cdot\U_+^k\sset \U_+^{k+n}\cdot A,\quad\forall k\gg 0\,.
\end{equation}
Recall further that for any left $A$-module $M$,
the linear dual $M^*=\Hom_{_\C}(M,\C)$ has a natural
right $A$-module structure.
Given a {\it left} $A$-module $M$ we define
\begin{equation}\label{Jdual}
\j^\dag(M)\,:=\,\lim_{\stackrel{\longrightarrow}{k}}\, \bigl(M/\U_+^k\cd
 M\bigr)^*=\{m^*\in M^*\;\mid\; m^*\bigl(\U_+^k\cdot M\bigr)=0\enspace
\text{for some}\enspace k=k(m^*)\gg 0\}.
\end{equation}
Formula \eqref{continuous} insures that $\j^\dag(M)$ is an
$A$-submodule of the right $A$-module $M^*$. Moreover,
it is clear that the action on $\j^\dag(M)$
of the Lie subalgebra $\n\sset A$ is locally nilpotent.
Thus, we get
a functor $\j^\dag:\,\text{\it left $A$-modules}\too
\text{\it right $A$-modules locally-nilpotent relative to $\n$}$.

It is convenient to get the following slightly different interpretation
of this functor. Observe that formula (\ref{continuous}) 
implies
 that the multiplication map: $ A\stackrel{a\cdot}{\too}A$
(by any element   $a\in A$)
is continuous in the topology on $A$ defined by
the set $\{\U_+^k\cdot A\}_{k\geq 1}$ of fundamental neighborhoods of
zero. Thus, the completion
$\hA=\lim\limits_{\stackrel{\longleftarrow}{k}}\, A/\U_+^k\cd A$
acquires the structure of a complete topological algebra.
There is a canonical algebra map:
$A\to \hA$ with dense image.
Similarly, for any left $A$-module $M$, the completion
\begin{equation}\label{Jhat}
{\widehat \j}(M)
= \lim_{\stackrel{\longleftarrow}{k}}\;
M/ \U_+^k\cd M
\end{equation}
acquires the structure of a complete topological 
left $\hA$-module. 
A noncommutative analogue of the
standard Artin--Rees lemma implies that, on the category of
finitely-generated
$A$-modules, the  functor
$M\mapsto {\widehat  \j}(M)$ is  exact.
Furthermore, if $M$ is finitely generated over
$A$, then  one has
a canonical surjective morphism  $\hA\otimes_{_A}M\onto
{\widehat  \j}(M)$, see \cite[ch.10]{AM}.
It is clear that $\j^\dag(M)$ is
nothing but the {\it continuous dual} of
${\widehat  \j}(M)$.

The  object $\widehat \j(M)$
is too
large to
be finitely generated over $A$. Using the element
$\bh\in A$, that has not played any role so far,
one  may do better provided one works with
classed of $A$-modules considered below.


\subsection{}
Let $\Amod$ be the category of left $A$-modules which
are finitely generated as $\U_+$-modules. This is
an abelian category since $\U_+$, being a quotient of
the augmentation ideal in $\U\n$, is Noetherian.
Notice that for any $  M\in \Amod$ and any $k\geq 1$, 
the space $M/ \U_+^k\cd M$ is finite  dimensional.
 The action of   $\delta\in A$ on $M$
induces a $\delta$-action on each finite-dimensional space in the
inverse system
\begin{equation}
\label{inv_sys}
M/\U_+\cd M\;
\longleftarrow\; M/\U_+^2\cd M \;   \longleftarrow
\;M/\U_+^3\cd M\; \longleftarrow \;\ldots
\end{equation}

One deduces from the positivity of $\;\ad\delta$-eigenvalues on ${\n}$
that, for each $\lambda \in \C$, the corresponding inverse system
of generalized 
$\la$-eigenspaces $\,\big[M/\U_+^k\cd M\big]_\lambda$
 {\sl stabilizes},
i.e.,
the projections in \eqref{inv_sys} induce isomorphisms
$\displaystyle \,\big[M/\U_+^{k+1}\cd M\big]_\lambda
\iso\big[M/\U_+^k\cd M\big]_\lambda,\,$ for all $k$ sufficiently large.
\begin{definition}\label{jacquet}
 Let
 $\j(M)=\bigoplus_{\la\in\C}\; \big[\lim\limits_{\stackrel{\longleftarrow}{k}}\;
M/ \U_+^k\cd M\big]_\lambda$ be the 
direct sum of all
such "stable" generalized $\delta$-eigenspaces,
to be called the Jacquet module for $M$.
\end {definition}
It is clear that $\j(M)$ may be identified with
the direct sum of all generalized 
$\delta$-eigenspaces in $\widehat \j(M)$. Thus, we have:
\smallskip

\noindent
$(1)\quad$ $\j(M)$ is an $A$-submodule in $\widehat \j(M)$ which
is dense in $\widehat \j(M)$ in the $\U_+$-adic topology.

\noindent
$(2)\quad$ $\delta$-action on $\j(M)$ is locally finite, with
finite  dimensional
generalized eigenspaces.

\noindent
$(3)\quad$  $\j(M)$ is finitely-generated over $\U_+$, that is,
$\j(M)\in \Amod$.

\noindent
$(4)\quad$ The set of $\delta$-eigenvalues in $\j(M)$ is
bounded from below by some constant $
\lambda(M) \in \C$.

\noindent
$(5)\quad$ $\j^\dag(M)=\bigl(\j(M))^\star$ is the  restricted dual of  $\j(M)$.
\smallskip

Properties (1) and  (2) above follow from   stabilization of
the eigenspaces in the inverse system (\ref{inv_sys}).
To prove (3), let $N_0\subsetneq N_1\subsetneq 
\ldots $ be a strictly  ascending chain of $\U_+$-submodules
in $\j(M)$.
Property (2) implies that, for each $i\ge 0$, there exists $\la(i)\in\C$
such that $\big[N_i\big]_{\la(i)} \neq \big[N_{i+1}\big]_{\la(i)}$.
It follows that for the corresponding closures in
$\widehat \j(M)$ we have 
$\overline{N}_i\neq \overline{N}_{i+1}$.
Thus, $\overline{N}_0\subsetneq 
\overline{N}_1\subsetneq 
\ldots $ is  a strictly ascending chain of $\widehat \U_+$-submodules
in $\widehat \j(M)$. But this contradicts the fact
that $\widehat \j(M)$ is a  Noetherian  $\widehat \U_+$-module
(the latter holds since
$M$ is finitely-generated
over $\U_+$, hence $\widehat \j(M)$ is finitely-generated
over $\widehat \U_+$). Thus, $\j(M)$ must be a finitely-generated
$\U_+$-module, and (3) is proved.
Property (4) follows from (3), and property (5) is clear from
(2) and (4). 
\smallskip

Since taking (generalized) 
eigenspaces is an exact functor, we conclude that the assignment
 $M\mapsto \j(M)$ gives an exact functor $\j:\ \Amod\to\Amod,$
called the Jacquet functor.

\subsection{} 
We say that $A$ has a triangular structure
if it  has either a commutative or a noncommutative
triangular structure.
From now one we fix a
triangular structure  $(A^\pm,\bh)$ on $A$.
Thus, one has categories  $\ou,\od$.
We also consider the category $\Amod$ for
the subalgebra $\U_+:=\bar A^+.$

\begin{lemma}\label{duality}
 \vi The  Jacquet functor gives
an exact functor $\j:\;\; \Amod\too \ou$. 

\vii
The functor $\j^\dag$ gives an
exact functor:
 $\Amod\too \od$ and, moreover, $\j^\dag(M)=\bigl(\j(M)\bigr)^\star,$
for any $M\in\Amod$.
\end{lemma}

\begin{proof} To prove (i) we note that, for any $M\in \Amod,$
the $A^-$-action on  $\j(M)$ is 
 locally-nilpotent since
the set of the real parts of $\delta$-eigenvalues in 
$\j(M)$ is
bounded from below by some constant $\lambda(M) \in \mathbb{R}$.
We know also that $\j(M)$ is 
 finitely-generated over $A^+$ (property (3) above), hence, over $A$.
Thus, $\j(M)\in\ou$,
and (i) is proved. 

Now, Lemma \ref{dual} yields
$\bigl(\j(M)\bigr)^\star\in \od$. But  property (5) above
says that
 $\j^\dag(M)\simeq\bigl(\j(M)\bigr)^\star.$
This
completes the proof of part (ii).~\end{proof}
\smallskip

Next, we introduce an abelian category $\AAA$  of 
left $A\otimes
A\opp$-modules
as follows.
First, we put
$\fu:=\bar A^+\otimes (A^-)\opp+ A^+\otimes (\bar A^-)\opp$,
a nonunital  subalgebra of $A\otimes A\opp$.
In the case where the  triangular
 structure is {\it commutative}, we let $\AAA$ be the category of
 $A\otimes A\opp$-modules
which are finitely generated over the subalgebra
$\fu$.

Now, let   the  triangular
 structure on $A$ be {\it noncommutative}.
Then the filtration on $A$ induces one on
$A\otimes A\opp$. Restricting the latter filtration
to the subalgebra $\fu
\sset A\otimes  A\opp$ makes $\fu$ a filtered algebra.
Taking associated graded algebras, one obtains
 a  graded algebra imbedding
$\gr \fu \into\gr(A\otimes A\opp) $.
We define  $\AAA$ to be the category of 
finitely-generated left $A\otimes A\opp$-modules $M$ 
such that $\gr M$ is a finitely
generated   $\gr\fu$-module,
for any {\it good} increasing
filtration on $M$. 
In particular, any object of $\AAA$ is a  finitely-generated 
$\fu$-module.

Recall  that, given a noncommutative triangular structure on $A$,
 we have the Lie subalgebras $\n^\pm\sset\bar A^\pm$.
In the case of a {\it commutative} triangular structure  on $A$,
we let $\n^\pm$ be any $\bh$-stable finite  dimensional
subspace of $\bar A^\pm$ that generates  $\bar A^\pm$ as an associative
algebra. We then regard  $\n^\pm$ as abelian  Lie subalgebras in $A$.
Thus, for  any  (commutative or  noncommutative) triangular structure on
$A$, we have the Lie algebra
 $\,\n^+\otimes 1 + 1\otimes \n^-\sset A\otimes A\opp$
such that $\fu$, an associative subalgebra of the algebra
 $A\otimes
A\opp$,
is a quotient of the augmentation ideal
of ${\mathcal U}(\n^+\oplus\n^-)$, the  enveloping algebra of the Lie algebra
$\n^+\oplus\n^-\cong
\n^+\otimes 1 + 1\otimes \n^-$.
Thus, on the category $\AAA$,
 there is a well  defined Jacquet
functor relative
to the Lie algebra $\n^+\otimes 1 + 1\otimes \n^-\sset  A\otimes A\opp$
and the element $\bh\otimes 1 + 1\otimes\bh$.
Furthermore, one mimics the proof of Lemma \ref{duality}
to show an analogue of the Lemma, saying that the Jacquet
functor is exact and takes the category 
$\AAA$ into $\oa$. Thus, we get an exact  functor 
\[\j:\, \AAA\too\oa.\]

\begin{proposition}
\label{faithful} 
\vi The Jacquet functor $\j:  \AAA\too\oa$
is {\sl faithful}, that is, we have 
$\,\,\dis M\in \AAA\,\,\text{and}\,\,
M \neq 0 \;\; \Longrightarrow \;\; \j(M) \neq 0.
$

\vii Any object of the category $\AAA$ has finite length.
\end{proposition}

\begin{remark} There is also a "one-sided" analogue of this proposition
formulated as follows.
Assume, for concreteness, that the  triangular 
structure on $A$ is  noncommutative. We
 define $\AM$ to be 
 the category of 
finitely-generated left $A$-modules $M$
such that $\gr M$ is a finitely
generated   $(\gr A^+)$-module
for any {\it good} increasing
filtration on $M$.
In particular, any object of $\AM$ is a  finitely-generated 
$\bar A^+$-module,
but the category  $\AM$ is, in general, {\it smaller}
than the category $\Amod$ for $\U_+=\bar A^+$. 

Similarly to  Proposition \ref{faithful}
one proves
\begin{proposition}
\label{faithful2}
\vi The Jacquet functor gives a {\sl faithful}
 functor $\j:  \AM\to\ou$.

\vii  Any object of the category $\AM$ has finite length.\qed
\end{proposition}
\end{remark}

\subsection{} Our proof of Proposition \ref{faithful} (and  Proposition \ref{faithful2})
will be based on  Gabber's Separation theorem, which we now
recall, cf. [Ga] for more details.

Let $\r=\n\rtimes\C\cd\bh$ be an arbitrary finite dimensional solvable Lie
algebra such that
the adjoint $\bh$-action on $\n$ is diagonal
with all the eigenvalues in $\Z_+$,
 as at the beginning of this section.
Write $\U_\bullet\r$ and $\U_\bullet\n$
for the standard increasing filtrations on the corresponding
enveloping algebras. Thus, $\gr(\U\n)= \sym{\n}$.

One has
the following result due to O. Gabber, [Ga], Theorem 1.

\begin{theorem}[Separation theorem]
\label{gabber}  
Let $M$ be a nonzero $\U{\mathfrak r}$-module, and
  $\{M_j\}_{j \geq 0}$ an 
increasing filtration on $M$ compatible with the
 $\U{\mathfrak r}$-action (i.e., such that
$\U_i\r\cdot M_j\sset M_{i+j}\,,\,\forall i,j\geq 0$).
Assume, in addition, that $\gr M$ is finitely generated {\sl
  over $\gr(\U{\n})$} (not only over $\gr(\U{\r})$).
Then, we have  $\,{\n} \cdot M \neq M$.
\end{theorem}

\begin{remark}
If  ${\n}$ is abelian, then   the separation theorem simplifies,
and becomes the following  standard 
result  in Commutative Algebra (cf. e.g. [SW]):
{\it If $M$ is an $\U\r$-module which is finitely-generated
over the subalgebra $\U\n$, then   $\,{\n} \cdot M \neq M$}.
(no filtration is needed in this case).
To prove this,
we claim first that the point $0
\in \Spec(\sym{\n})$ belongs to $\Supp(M)$. If not, then there is a
polynomial $P \in \sym{\n}$ that vanishes on $\Supp\ M$ and such that
$P(0) \neq 0$. Replacing $P$ by its high enough
power we may achieve that $P$
annihilates $M$, i.e., $P \in \ann(M) \subset \sym{\n}$. But the space
$\ann(M)$ is clearly stable under the adjoint $\delta$-action on
$\sym{\n}$. Moreover, since $P = P(0) + P_1$, where $P_1 \in
 {\n}\!\cdot\!(\sym{\n})$, and all weights of $ad \delta$-action 
on ${\n}\!\cdot\!\sym{\n}$ are
strictly positive, we deduce from $P \in \ann(M)$ that $P(0)$, $P_1 \in
\ann(M)$. Since $P(0) \neq 0$ this yields $1 \in \ann(M)$, a
contradiction. Thus we have proved $0 \in \Supp \ M$.
 But then $M/\n\cd M$, the geometric fiber of $M$ at $0$,
is nonzero, due to the Nakayama lemma.
Thus, $M\neq \n\cd M$. $\enspace\Box$
\end{remark}

Here is an example\footnote{It was kindly communicated to me
by O. Gabber} showing that
the separation result may fail, in general, if (for $\n$ noncommutative)  the
assumption of  Theorem \ref{gabber}:   "$\gr M$ {\it is finitely generated
  over} $\gr(\U{\n})$" is replaced by the
weaker assumption:  "$M$  {\it is finitely generated
  over} $\U{\n}$".
\smallskip

\noindent
{\bf Example.\;}
Let $\n$ be the 3-dimensional Heisenberg Lie algebra with basis
$x,y,z,$ where $z$ is central and $[x,y]=z$. Define
a semidirect product $\r=\n\rtimes\C\cd\bh$
by the commutation relations $[\bh,x]=x\,,\,[\bh,y]=y,$
and $[\bh, z]=2z$. 
Further, let $D$ be the associative algebra
 of polynomial differential operators on the 2-plane
$\BA:=\C^2,$ with coordinates $y,z$.
The assignment  
$$x\mapsto z\frac{\partial}{\partial y}\enspace,
\enspace y\mapsto y\enspace,
\enspace z\mapsto z\enspace,\enspace\bh\mapsto 2z\frac{\partial}{\partial z}
+y\frac{\partial}{\partial y}\,.
$$
extends to a Lie algebra imbedding: $\r\into D.$ 
We set $f(y,z):= y\cdot z -1$. Thus $f\in \C[y,z]$ is
 a polynomial, and
 the equation $f=0$ defines a hyperbola
 $X\sset
\BA$. It is clear that

\noindent
$\bullet\quad$ The curve $X$ is smooth, irreducible, and does not intersect the axis $z=0$, and

\noindent
$\bullet\quad$ The restriction to $X$ of the projection $\BA\onto \C$
along the $y$-axis gives an \'etale morphism
$X\to \C$.

Set $U:=\BA\ssminus X$, an affine Zariski open subset in $\C^2$.
The vector space $M:=\C[U]\big/\C[\BA]= \C[\BA][\frac{1}{f}]\big/\C[\BA]$ has
a natural $D$-module structure (that makes $M$ a simple holonomic
$\dd$-module on the plane $\BA$, with support $X$). It is straightforward to verify that
the action of $\frac{\partial}{\partial y}$ on $M$
induces isomorphisms:
$$\mbox{$\frac{1}{f^n}$}\cdot\C[y,z]\Big/\mbox{$\frac{1}{f^{n-1}}$}\cdot\C[y,z] \iso
\mbox{$\frac{1}{f^{n+1}}$}\cdot\C[y,z]\Big/\mbox{$\frac{1}{f^n}$}\cdot\C[y,z]
,\quad\forall n=1,2,\ldots\,.
$$
It follows easily that $M=\U\n\cdot[\frac{1}{f}],$ i.e.,
the class $[\frac{1}{f}]\in
\C[\BA][\frac{1}{f}]\big/\C[\BA]$ 
generates $M$ over the subalgebra $\U\n\sset D$.

On the other hand, it is clear that $z\cdot M=M$,
i.e., the separation property {\sl fails} for $M$
(in this case,  $\gr(M)$ is not
finitely generated over $\gr(\U\n)$). $\lozenge$

\subsection{Proof of Proposition \ref{faithful}.} The argument is quite standard.
Set $\n=\n^+\otimes 1\oplus1\otimes\n^-$, a Lie subalgebra
in $A\otimes A\opp$. Thus, the algebra $\fu$  is a quotient of the augmentation ideal
of the enveloping algebra of this Lie algebra.
We have the Jacquet functor
$\widehat \j$ on $A\otimes A\opp$-modules
defined as  the completion with respect to the
subalgebra $\fu$.
In view of the identifications above, 
for any $M\in \AAA$, we can write:
${\widehat \j(M)}= \lim\limits_{\stackrel{\longleftarrow}{k}}\;
M/ \fu^k\cd M$.
The kernel of the canonical
map $M\to \widehat \j(M)$ is clearly equal to
$K=\bigcap_{k\geq 0}\, \fu^k M$.
It is routine to verify
 that $K$ 
 is an $A\otimes A\opp$-submodule in $M$.
Moreover, a noncommutative version of the Artin--Rees
lemma, see [AM], [Ca], implies that
$\fu\cdot K=K$ or, equivalently,
that $\n\cdot K=K$.

We would like to apply the separation 
theorem to deduce that
$K=0$. Observe that $K \in \AAA$,
since $M\in \AAA$.
In the case where the triangular structure $(A^\pm,\bh)$ is commutative,
it follows readily that all the conditions on  $K$
required by  the Separation Theorem  \ref{gabber} 
hold trivially. In the case of a noncommutative
 triangular structure $(A^\pm,\bh)$, choose a good
filtration $K_\bullet$
on $K$. Then, since $\bh\in A_1,$ for any $i\in \Z,$ we
get  
$
(\bh\otimes 1 + 1\otimes\bh)\cdot K_i\sset K_{i+1}.$
Further,  Lemma \ref{nc} guarantees that $\gr K$
is a finitely-generated $\gr\fu$-module.
Thus, all the conditions of the Separation Theorem  \ref{gabber} 
hold. Thanks to the theorem,
we conclude that $K=0$.

It follows from the above that the canonical map $M \to 
{\widehat \j(M)}$ is injective.
 Hence, for any nonzero $M\in \AAA$
we have ${\widehat \j(M)}\neq 0$.
We conclude that $\j(M)$, being a dense subspace in 
${\widehat \j(M)}$, is also nonzero,
and part (i) follows.

To prove (ii), it suffices to  verify 
the descending chain condition for
any object of the category
$\AAA$. Let
$M=M^0\supsetneq M^1\supsetneq  M^2\supsetneq  \ldots$ be a descending chain
of submodules in $M$.
Using  the exactness of the
Jacquet functor and part (i), from
 $M^i/M^{i+1}\neq 0\,,\,
\forall i,$ we deduce that
$\j(M^i)/\j(M^{i+1})=\j(M^i/M^{i+1})\neq 0$.
Hence,
$\j(M^0)\supsetneq  \j(M^1)\supsetneq  \j(M^2)\supsetneq   \ldots$ is a 
strictly  descending chain
of submodules in $\j(M)$. 
This contradicts the fact
that $\j(M)$, being an object of $\oa$,
must have finite length by Lemma
\ref{triv}.
Claim (ii) follows.~\qed

\subsection{Proof of Theorem \ref{duflo}.} First of all, we observe
that conditions $(T0)$, $(T0')$ insure that $A$, the diagonal
 $(A\td A)$-bimodule, is an object of category $\AAA$.
Thus, the last statement of Theorem \ref{duflo}  follows
directly from Proposition \ref{faithful}(ii).

Now, we turn to other  statements of the theorem. The implication
(i)$\,\Longrightarrow\,$(ii) is standard.

To prove  (ii)$\,\Longrightarrow\,$(iii),
fix a prime ideal $I\sset A$. As we have observed above,
 $A\in \AAA$, hence  $A/I\in \AAA$. Therefore,
 we may apply the Jacquet functor
to get an object $\j(A/I)\in \oa$.

Given an $A\otimes A\opp$-module, resp. $(A\td A)$-bimodule,
$M$, we write $\lann(M)$ for the annihilator of $M$ in the
subalgebra $A\otimes 1 \sset A\otimes A\opp$.
The proof of  Proposition \ref{faithful}
shows that $A/I$ maps injectively into $\widehat{\j}(A/I)$
and, moreover, the image of $A/I$ is dense in $\widehat{\j}(A/I)$.
Thus we  have
$I=\lann(A/I)=\lann\bigl(\widehat{\j}(A/I)\bigr).$
Furthermore, $\lann\bigl(\widehat{\j}(A/I)\bigr)=
\lann\bigl({\j}(A/I)\bigr),$ for $\j(A/I)$ is dense
in $\widehat{\j}(A/I)$. Thus,
$\lann\bigl({\j}(A/I)\bigr)=I$ is a prime 
ideal in $A$.\footnote{It was pointed out by a referee that
a similar observation is contained in [J2, n.2.4].}

Now,  any object of $\oa$, in particular,
$M=\j(A/I),$  has finite length by Lemma \ref{triv}.
Let $0=M_0\sset M_1\sset M_2\sset\ldots\sset M_n=M$
be its Jordan--H\"older series.
For each $i=1,2,\ldots,n,$ put $J_i:=\lann(M_i/M_{i-1}).$
Then,  one proves by a standard argument that
there exists  $m\in [1,n]$ such that one has $I=J_m$, see e.g. [BG].
Below, we provide more details for the reader's convenience.

First, since $J_i(M_i)\sset M_{i-1}$ we deduce that
$$
J_1\cdot\ldots\cdot J_{n-1}\cdot J_n(M)\;\;\sset\;\;
J_1\cdot\ldots\cdot J_{n-1}(M_{n-1})\;\;\sset\;\;\ldots\;\;\sset\;\;
J_1(M_1)\;=\;0.
$$
Hence, the ideal $J_1\cdot\ldots\cdot J_n$ annihilates $M$.
Therefore, we have $J_1\cdot\ldots\cdot J_n\sset \lann(M)$.
But, $\lann(M)=I$ being prime, the inclusion above implies
that there exists an $m$ such that $J_m\sset I$.
On the other hand,
 since $M_i/M_{i-1}$ is a subquotient of $M$,
 for any $i=1,2,\ldots,n,$ one clearly has an opposite
inclusion $I=\lann(M)\sset\lann(M_i/M_{i-1})=J_i$.
Thus, we obtain $I=J_m=\lann(M_m/M_{m-1})$, as claimed.

By construction,
 $M_m/M_{m-1}$ is a simple object
of $\oa$. Hence, Lemma \ref{triv} yields\linebreak
$M_m/M_{m-1}= L'\boxtimes L'',$
where $L'\in \ou$ and $L''\in \od$ are some simple objects.
Further, it is clear that $\lann(L'\boxtimes L'')=
\lann(L')$. Thus, we deduce that
$\,I=\lann(M_m/M_{m-1})=$
$\lann(L'\boxtimes L'')=\lann(L')\,$ is the annihilator
of a simple object of $\ou$. Part (iii) follows.
The implication
(iii)$\,\Longrightarrow\,$(i) is trivial.
\qed

\section{Applications to symplectic reflection algebras}
\subsection{} 
Let $(V,\om)$ be a finite-dimensional
symplectic vector space, and $\G\sset Sp(V,\om)$
a finite subgroup of symplectic automorphisms of $V$.
An element $s\in \G$ is called  a {\sl symplectic
reflection} if $\rk(\id-s)=2.$
Let $S$ denote the set of symplectic reflections in $\G.$
The group $\G$ acts on $S$ by conjugation.
For each $s \in S$, there is an $\om$-orthogonal direct sum decomposition
$V= {\mathtt{Image}}(\id-s) \oplus \Ker(\id-s)$,
and we write $\om_s$ for the (possibly degenerate)
skew-symmetric form on $V$ which coincides with $\om$
on ${\mathtt{Image}}(\id-s)$, and has $\Ker(\id-s)$ as its
 radical.

Write $\C\G$ for
the group algebra of $\G$, and  $(TV)\#\G$ 
 for the cross product of the tensor algebra $TV$ with
$\C\G$. From now on, we identify $V$ with $V^*$ via the symplectic form.
Thus, the symmetric algebra on $V$ is identified with
 $\C[V]$, the polynomial algebra on $V$.

In [EG], the authors have introduced a class of associative algebras
 $\sH_{t,c}(V,\om,\G)$, called symplectic
reflection algebras, as follows.
Fix an ${\mathtt{Ad}^{\,}}\G$-invariant function $c: S \to \C\,,\, s\mapsto
c_s,$ and a complex number $t\in\C$.
Let  $\kappa: V\times V \to \C\G$ be a
skew-symmetric $\C$-bilinear pairing given 
by the formula
$$\kappa(x,y) = t \cdot \om(x,y)\cdot 1 + \sum\nolimits_{s\in S}\; c_s\cdot
\om_s(x,y)\cdot s,\quad\forall x,y \in V\,.$$
We define the symplectic
reflection algebra with parameters $(t,c)$ by
\begin{equation}\label{rel}
\sH_{t,c}(V,\om,\G)
\; :=\; (TV\#\G)/I\langle x\otimes y - y\otimes x - \kappa(x,y) \in
T^2V\, \oplus\, \C\G\rangle_{\,x,y \in V}\,\,,
\end{equation}
where 
 $I\langle\ldots\rangle $ stands for the two-sided ideal in
$(TV)\#\G$ generated by the indicated set.
 Thus, $\sH_{t,c}(V,\om,\G)$
is an associative algebra which may be thought of as a
deformation of $\C[V]\#\G$ (= cross-product of $\C[V]$
with $\G$). Clearly,  $\sH_{t,c}(V,\om,\G)$ contains $\C\G$
  as a subalgebra. Furthermore,
there is a natural increasing filtration on
 $\sH_{t,c}(V,\om,\G)$ such that $\C\G$ has filtration degree zero,
and elements of $V\sset TV$ are assigned  filtration degree one.
It has been proved in  [EG] that, for
the corresponding associated graded  algebra, there is 
a canonical  graded  algebra isomorphism
\begin{equation}\label{H_filt}
\gr\sH_{t,c}(V,\om,\G)\;\simeq\;\C[V]\#\G \qquad\text{({
\sl{Poincar\'e--Birkhoff--Witt property}})}.
\end{equation}

Write $\e=\frac{1}{\sharp\G}\sum_{g\in\G}\,g\in \C\G$
for the symmetrizer idempotent, viewed as an element in
$\sH_{t,c}(V,\om,\G)$.
We let ${\sf{eHe}}:=\e\cdot\sH_{t,c}(V,\om,\G)\cdot\e$
be the {\it spherical subalgebra} in  $\sH_{t,c}(V,\om,\G)$.
The increasing filtration on
 $\sH_{t,c}(V,\om,\G)$ induces a filtration
on the  spherical subalgebra. From (\ref{H_filt}), one finds that
$\gr({\sf{eHe}})=\C[V]^\G=\C[V/\G],$ the algebra
of $\G$-invariant polynomials on $V$. 
The canonical  Poisson structure on the commutative
algebra $\gr({\sf{eHe}})$ makes $V/\G$ a  Poisson variety.
It is
easy to see, cf. ([EG], Theorem 1.6),
 that 
this Poisson structure on $V/\G$
equals
the $t$-multiple of  the standard one induced
by 
the symplectic structure on $V$
(for $t=0$ the spherical subalgebra ${\sf{eHe}}$ becomes commutative,
and the corresponding Poisson structure on
$\gr({\sf{eHe}})=\C[V/\G]$ reduces to zero).

\subsection{Symplectic leaves in $V/\G$.}
The  leaves
of the standard Poisson structure on $\C[V/\G]$ 
are     described
as follows, see  [BrGo]. Let $\Isot$ denote the
(finite) set of all subgroups $G\sset \G$ that occur as 
isotropy groups of points in $V$.
Given a  subgroup $G\sset \G$,
write $V^G\sset V$ for the vector subspace of
$G$-fixed points. Let 
$V_G$ denote the
image of $V^G$ under the projection $V\onto V/\G$,
and let $\overset{_\circ}{V}_G= V_G\ssminus\bigcup'_{_{G'\in\Isot}}\; 
\left(V_G\cap V_{G'}\right)$
be the complement in $V_G$ of the union of
all {\it proper} subsets in  $V_G$ of the form $V_G\cap V_{G'}$
(here 
the symbol $\cup'$ indicates that the union is taken
over proper subsets only). 
It is easy to see that 
for each subgroup  $G\in \Isot$, the set
$\overset{_\circ}{V}{}^{_G}$ equals
the image of the
subset $\overset{_\circ}{V}{}^{_G}\sset 
V^G$ formed by the points whose isotropy group equals
$G$. If $N_G$ denotes the normalizer of $G$ in $\G$,
then the group $N_G/G$ acts freely on $\overset{_\circ}{V}{}^{_G}$
and we have $\overset{_\circ}{V}_G\simeq
\overset{_\circ}{V}{}^{_G}\big/(N_G/G)$.
It follows that each set $\overset{_\circ}{V}_G$ is
a smooth connected locally-closed
subvariety of $V/\G,$ and one shows that these varieties
 are exactly the  symplectic leaves in
$V/\G$, see e.g., [BrGo]. Hence, there are only finitely
many  symplectic leaves. 

    From the Irreducibility Theorem \ref{irreducibility}
we deduce

\begin{corollary}\label{symplectic_refl}
Let   $\sH_{t,c}(V,\om,\G)$ be a  symplectic
reflection algebra, and $t\neq 0$.
Then, for any primitive ideal $I$ in
the spherical subalgebra $\e\cd\sH_{t,c}(V,\om,\G)\cd\e$,
the variety $\var(I)$ has the form $V_G$, for a certain
subgroup $G\in\Isot$.
\end{corollary}

We remark that Theorem \ref{irreducibility}
is not applicable to the symplectic reflection
algebra $\sH_{t,c}(V,\om,\G)$ itself since
$\gr\sH_{t,c}(V,\om,\G)=
\C[V]\#\G,$ is a noncommutative algebra.
Note however, that 
finitely generated $\C[V]\#\G$-modules may be naturally
identified with $\G$-equivariant coherent sheaves on
$V$. The support of such a sheaf is a  $\G$-stable
subvariety in $V$.
In particular, given a two-sided ideal $I\sset \sH_{t,c}(V,\om,\G)$
one can view $\gr\bigl(\sH_{t,c}(V,\om,\G)/I\bigr)$
as a finitely generated $\C[V]\#\G$-module, via (\ref{H_filt}).
Let $\var(I)\sset V$ stand for its support.

Using the technique of {\it Poisson orders} developed in [BrGo],
one can refine our argument to obtain the following result.

\begin{proposition}\label{symplectic_refl2}
For  any primitive ideal $I\sset\sH_{t,c}(V,\om,\G)$,
the variety $\var(I)$ is the $\G$-saturation
of vector subspace $V^G$, for a certain
subgroup $G\in\Isot$.\hfill\qed
\end{proposition}

\subsection{Rational Cherednik algebras.}
Let $ W $ be a finite Coxeter group in a
complex vector space $\h$.
Thus, $\h$ is the complexification
of a  real Euclidean  vector space,
with inner product $(-,-)$, and  the group $W$ is
generated by a finite set $S\subset W$ of
reflections $s\in S$ with respect to certain
hyperplanes $\{H_s\}_{s\in S}$
in that  Euclidean  space.

For each $s\in S$, we choose 
a nonzero
linear function $\alpha_{s}\in\h^*$ that vanishes on $H_s$
(called the {\it positive root} corresponding to $s$),
and let $\alpha_s^{\vee}=2(\alpha_s, -)/(\alpha_s,\alpha_s) \in \h $
be the corresponding coroot.
The group $W$ acts   on the set $S$,
and  also acts diagonally  on $V=\h\oplus\h^*$,
 by conjugation. We  equip  the space
 $V=\h\oplus\h^*$  with the canonical symplectic
structure.

 Following [EG], to each 
 $W$-invariant function $c: S \to \C\,,\, c\mapsto c_s,$
 one associates the symplecic reflection algebra $\sH_c:=\sH_{1,c}(V,\om,W)$,
 called  the {\it rational Cherednik algebra}.
The algebra $\sH_c$
is generated
by the vector spaces $\h$, $\h^*,$ and the set
$W,$ with defining relations  
(cf.  formula (1.15) of [EG] for $t=1$) given
by
\begin{equation}
\label{defrel}
\begin{array}{lll}\displaystyle
&{}_{_{\vphantom{x}}}w\cd x\cd w^{-1}= w(x)\;\;,\;\;
w\cd y\cd w^{-1}= w(y)\,,&\;
\forall y\in \h\,,\,x\in \h^*\,,\,w\in W\break\medskip\\
&{}^{^{\vphantom{x}}}{}_{_{\vphantom{x}}}x_1\cd x_2 = 
 x_2\cd x_1\;\enspace,\;\enspace
y_1\cd y_2=y_2\cd y_1\,, &\;
\forall y_1,y_2\in \h,\;x_1\,,\,x_2 \in \h^*\,\break\medskip\\
&{}^{^{\vphantom{x}}}y\cd x-x\cd y = \langle y,x\rangle
+\sum\limits_{^{_{s \in S}}}\,
c_s\cd\langle y,\alpha_s\rangle
\langle\alpha_s^\vee,x\rangle \cd s\,,&\; \forall y\in
\h\,,\,x\in \h^*\,.
\end{array}
\end{equation}

 Thus, the elements  $x\in\h^*$ generate a 
subalgebra $\C[\h]\subset \sH_c$, of polynomial functions on $\h$,
 the elements  $y\in\h$ generate a 
subalgebra $\C[\h^*]\subset \sH_c$, of polynomial functions on 
the dual space,
and the elements $w\in W$
span  a copy of the group
algebra $\C W$ sitting naturally inside $\sH_c$.
Furthermore, the Poincar\'e--Birkhoff--Witt property yields
  a "triangular 
decomposition" $\,\sH_c\simeq \C[\h^*]\otimes\C[W]\otimes\C[\h],$
see ([EG], Corollary 4.4). We put $A^+:=
\C[\h]^W$ and
$A^-:=
\C[\h^*]^W$. Further, in  [BEG] we have constructed
an element $\bh\in \sH_c$ such that the data $(A^\pm,\bh)$ 
gives a commutative triangular structure on the
rational Cherednik algebra $\sH_c$.

Associated to the above data, we have
 categories
$\ou, \od$ (see Definition \ref{oo}).
Applying Theorem \ref{duflo} (note that we are using
here only the easy, purely commutative,
 part of Gabber's separation theorem, since
$A^\pm$ are commutative algebras), we obtain the following result

\begin{corollary}\label{cher} 
Duflo theorem
holds for the rational Cherednik algebra $\sH_c$.
\end{corollary}

In [BEG], the authors  introduced another category, $\O$, defined as the
category of finitely-generated left $\sH_c$-modules $M$
such that the $A^-$-action on $M$ is locally-finite.
It is clear from the triangular
decomposition $\,\sH_c\simeq \C[\h^*]\otimes\C[W]\otimes\C[\h]$
that, in the notation of \S4, we have 
$\ou\sset \O\sset \Amm$. Note that the inclusion
$\ou\sset \O$ is {\it strict} since
we do not require the  locally-finite
$\bar A^-$-action to be  locally-nilpotent.
We have also introduced in [BEG] a category $\HC$ of Harish-Chandra
bimodules over the Cherednik algebra $\sH_c$,
and showed that (in our present notation)
any object of $\HC$ belongs to $\AAA$.
Thus, from Propositions \ref{faithful}--\ref{faithful2}
we deduce
\begin{corollary}
\label{cher2}
Any object of the categories $\O$ and $\HC$ has finite length.
\hfill\qed
\end{corollary}

\begin{remark} In [BEG] we were only able to prove that
any object  of the category $\ou$ has finite length, but the
technique in loc. cit. was insufficient to
prove Corollary \ref{cher2}
in full generality.
\end{remark}

As has been explained before Proposition \ref{symplectic_refl2},
for any $\sH_c$-module, hence for an object $M\in \O$,
one
can define a $W$-stable subvariety
$\gsup{M}\sset \h\oplus\h^*$.
The  $A^-$-action on $M$  being locally-finite,
it follows, since $A^-= \C[\h^*]^W,$
that $\gsup{M}\sset \h\oplus \{0\}$.
We identify $\h^*$ with $\h$ via the invariant form and observe
that, since $W$ is finite,  the isotropy group  of a point $(y,x)\in  \h\oplus \h^*=
\h\oplus \h$ coincides with  the isotropy group of the
generic linear combination $t_1\cdot y + t_2\cdot x\,,\, t_1,t_2\in \C$.
Therefore we see that, in the notation of 
Theorem \ref{symplectic_refl}, 
the set $\Isotw$ is exactly the set ${\mathtt{Parab}}(W)$ of all parabolic
subgroups in $W$, cf. [Hu].
\begin{theorem}
\label{levi}
For  any simple object $M\in \O$,
the variety $\gsup{M}$ is the $W$-saturation
of vector subspace $\h^G\oplus \{0\}$, for a certain
subgroup $G\in {\mathtt{Parab}}(W)$.
\end{theorem}

\noindent
{\sl Sketch of Proof.\;} We adapt the known argument used
in the case of highest weight modules over $\U\g$,
the enveloping algebra of a semisimple Lie algebra.
To this end, write $\ann M\sset \sH_c$ for the annihilator
of $M\in \O$. Thus, $\ann M$ is a primitive ideal in $\sH_c$.

Introduce the notation $\Lambda:=\gsup M\sset \h\oplus\h^*.$
Clearly, $\La\sset \gsup(\sH_c/\ann M )$. 
Further, by Proposition \ref{symplectic_refl2} 
there exists a subgroup $G\in {\mathtt{Parab}}(W)$ such that
one has 
$\gsup(\sH_c/\ann M )=W\cdot (\h\oplus\h^*)^G$.
Therefore we get  
\begin{equation}\label{pf1}
\La\sset (\h\oplus\{0\})\;\cap\;W\cdot (\h\oplus\h^*)^G
=W\cdot (\h^G\oplus\{0\}).
\end{equation}

Next, we claim that 
\begin{equation}\label{pf2}
\dim \La =\dim (\h^G)\,.
\end{equation}
To prove this, choose a finite  dimensional
$\C[\h^*]^W$-stable subspace $M^0\sset M$ that generates
$M$ over the subalgebra $\C[\h]^W$. For each
$i\geq 0$, set $M_i= F_i\sH_c\cdot M^0$, where $\{F_i\sH_c\}_{i\geq 0}$
is the standard increasing filtration on the Cherednik algebra $\sH_c$,
cf. [EG]. This gives a filtration on $M$ compatible with
that on $\sH_c$ and such that $\gr M$ is
finitely-generated over the subalgebra $\gr\C[\h]^W \sset
\gr{\sH_c}.$ Hence, the restriction to $\gsup M\sset \Spec(\gr{\sH_c})$
of the canonical projection
$\Spec(\gr{\sH_c})\onto \Spec(\C[\h]^W)$ is a finite morphism.
Therefore the dimension of the variety $\gsup M$ is unaffected
by replacing the algebra $\sH_c$ by the algebra $\C[\h]^W$.
Now, following Joseph [J3], we define
an increasing chain of $\C[\h]^W$-submodules
 $D_k(M)\sset \End_{_\C}M\,,\,
k=0,1,\ldots,$ inductively as follows:
$$ D_0(M)=\C[\h]^W\,,
\quad\text{and}\quad
D_k(M)=\{u \in \End_{_\C}M\;\;\big|\;\;
[u\,,\,\C[\h]^W]\sset D_{k-1}(M)\}\,.
$$
Set $D(M):= \bigcup_{k\geq 0}\,D_k(M).$
Joseph showed that $D(M)$ is an associative subalgebra
in $\End_{_\C}M$ (which should be thought of as
an algebra of  \,"differential operators" on $M$).
Moreover, he proved  in ([J3], Lemma 2.3)
that $\dim\gsup D(M)\,\leq\, 2\dim(\gsup M)$.

Observe further  that the
 action in $M$ of an element $a\in \sH_c$ gives
an endomorphism $\hat{a}\in \End_{_\C}M$,
and the assignment $a\mapsto \hat{a}$ gives an 
algebra imbedding $\sH_c/\ann M \into \End_{_\C}M$.
The adjoint action of the subalgebra $\C[\h]^W$
on $\sH_c$ being locally-nilpotent, see [BEG],
the image of the embedding $\sH_c/\ann M \into  \End_{_\C}M$
is contained in $D(M)$.
      From this, following Joseph [J3], one derives
\begin{equation}\label{joseph}
\dim\gsup(\sH_c/\ann M )\;\leq\; \dim\gsup D(M)\; \leq\; 2\dim(\gsup M)\,.
\end{equation}
Hence, using \eqref{pf1} and writing "$\dim(-)$" instead of
"$\dim\gsup(-)$", we deduce
$$
\dim (\h\oplus\h^*)^G \;=\;
\dim(\sH_c/\ann M )\; \leq\; 2\dim M\;
=\; 2\dim\La\; \leq\; 2\dim\h^G\; =\;\dim (\h\oplus\h^*)^G \,.
$$
It follows that
 all the inequalities above must be equalities, and
\eqref{pf2} is proved. In particular, there exists 
 $\La^\circ\sset (\h^G\oplus\{0\})$, an irreducible
component of $\La,$ such that $\dim 
\La^\circ =\dim(\h^G\oplus\{0\})$. This dimension equality yields
$\La^\circ =\h^G\oplus\{0\}.$

We claim next that {\it all} irreducible components
of $\La= \gsup M$ have the same dimension, i.e.,
the following version of
 Gabber's equidimensionality theorem holds for the
Cherednik algebra.
\begin{proposition}\label{equi}
For any simple $\sH_c$-module $M$, the variety
$\gsup M$ is equidimensional.
\end{proposition}

To prove this, we recall that 
the equidimensionality theorem is, as explained e.g. in 
[Gi], a formal
consequence of  the  Gabber--Kashiwara theorem. 
The proof of the latter  theorem given,
e.g. in \cite[p.342--345]{Gi}, works for any 
filtered
algebra $\A$ such that $\gr\A$ is the coordinate ring
of a {\it smooth} affine algebraic variety.
Now, in our present situation,
we have  $\gr \sH_c = \C[\h\oplus\h^*]\#W$,
which is not a  commutative algebra.
However,
the formal microlocalization
construction  can be carried out with respect to any
multiplicative set $S\sset \A\,,\, (0\not\in S),$ for
any algebra $\A,$ provided
 the principal symbols of the elements of $S$
belong to the center of $\gr\A$,
see ([Gi], footnote on p. 337).
In our case we have a central subalgebra
$\C[\h\oplus\h^*]^W \sset \gr \sH_c.$
Furthermore, the  algebra $\gr \sH_c$
 has {\it finite} homological dimension,
as a cross-product of  the polynomial
algebra $\C[\h\oplus\h^*]$ with a finite group.
Going through the proof of the  Gabber--Kashiwara theorem
given in ([Gi], pp. 342--345) one sees that 
 the finiteness of homological dimension of $\sH_c,$
plus the existence of
the formal microlocalization with  respect to
 elements whose symbols belong to
$\C[\h\oplus\h^*]^W\ssminus\{0\}$, is enough
to conclude 
that the proof of   the   Gabber--Kashiwara theorem
applies to the algebra $\sH_c$
as well.
This implies Proposition \ref{equi}.

      From 
Proposition  \ref{equi} we deduce
 that all irreducible components
of the variety $\La=\gsup M$ have the same dimension, which is
equal to $\dim \h^G$. But then the inclusion in \eqref{pf1}
forces each irreducible component to be 
a $W$-translate of $\h^G$, and  Theorem \ref{levi} follows.
\qed
\medskip

\begin{remark}
It was pointed out to us by A. Yekutieli that
Gabber's equidimensionality theorem also
holds for the
spherical subalgebra  ${\sf{eHe}}$ in any symplectic reflection
algebra
${\sf{H}}=\sH_{t,c}(V,\om,\G)$, i.e., one
has\end{remark}
\begin{proposition} Let $\sH_{t,c}(V,\om,\G)$
be a  symplectic reflection
algebra. Then,
for any simple ${\sf{eHe}}$-module $M$, the variety
$\gsup M$ is equidimensional.
\end{proposition}

\noindent
{\sl Sketch of proof.\;} Recall that  PBW-proprty \eqref{H_filt}
implies that
$\gr({\sf{eHe}})=\C[V]^\G=\C[V/\G].$ We see that
 $\Spec\bigl(\gr({\sf{eHe}})\bigr) \simeq V/\G$ is not a smooth variety.
 Nonetheless, $ V/\G$
 is a Gorenstein variety, and
homological duality formalism
is known  to  work for Gorenstein varieties
as nicely as for smooth varieties. In particular, it follows from  [YZ]
that
the bimodule ${\sf{eHe}}$ is (up to shift) a  rigid dualizing complex
for the algebra  ${\sf{eHe}}$
(in the sense of [YZ]). Therefore, the proof
of   the   Gabber--Kashiwara theorem still goes through, cf. \cite[Purity~theorem]{YZ}.
\qed


\end{document}